\documentclass{article}
\usepackage{amssymb}
\usepackage{amsfonts}
\usepackage{amsmath}
\usepackage{graphicx}
\usepackage{hyperref}
\usepackage{authblk}

\newtheorem{theorem}{Theorem}

\newtheorem{definition}[theorem]{Definition}

\newtheorem{remark}[theorem]{Remark}

\begin{document}

\title{Analytical closed-form solution of the Bagley-Torvik equation}

\author[$\dagger$]{Juan Luis Gonz\'{a}lez-Santander}
\author[$\ddag$]{Alexander Apelblat}

\affil[$\dagger$]{Department of Mathematics, Universidad de Oviedo, 33007 Oviedo, Asturias, Spain.}
\affil[$\ddag$]{Department of Chemical Engineering, Ben Gurion University of the Neguev, 84105 Beer Sheva, Israel.}

\maketitle

\begin{abstract}
We calculate the solution of the Bagley-Torvik equation for arbitrary
initial conditions and arbitrary external force as the sum of two terms.
The first one is a linear combination of exponentials with error functions,
and the second one is a convolution integral
whose kernel is a linear combination of exponentials with error functions. The
derivation of the solution is carried out by using the Laplace transform
method and the calculation of a new inverse Laplace transform. The
aforementioned convolution integral can be calculated for the cases of a
sinusoidal- or a potential-type external force. In addition, we calculate the asymptotic behaviour of
the solution for $t\rightarrow 0^{+}$ and $t\rightarrow +\infty$.
The computation of this new
analytical solution is much faster and stable than other analytical
solutions found in the literature.
\end{abstract}

\section{Introduction}

\begin{figure}[htbp]
\centering\includegraphics[width=0.5\textwidth]{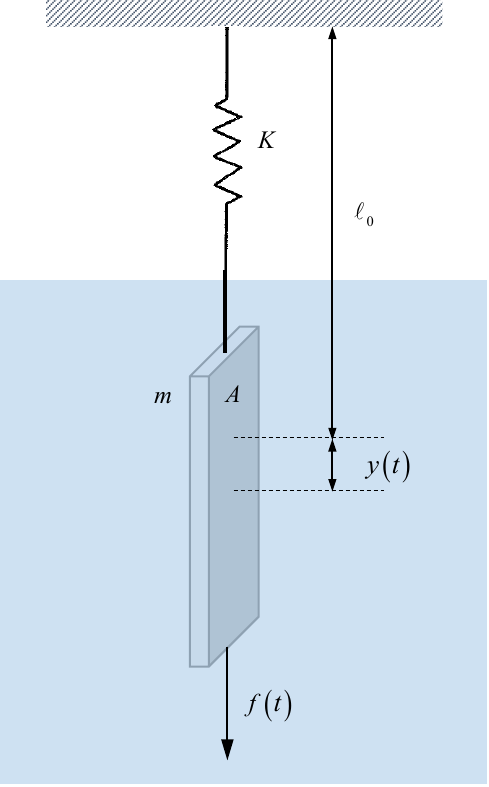}
\caption{Immersed plate in a Newtonian fluid connected by a massless spring.
}
\label{Figure: Bangley-Torvik setup}
\end{figure}

In 1984, P. J. Torvik and R. L. Bagley \cite{BangleyTorvik1984} considered
the motion of a rigid plate $y\left( t\right) $ of mass $m$ and area $A$
immersed in a Newtonian fluid of infinite extent with density $\rho $ and
viscosity $\mu $. This plate is connected by a massless spring of stiffness $%
K$ to a fixed point, where an external force $f\left( t\right) $ is applied
to it for $t>0$, see Figure \ref{Figure: Bangley-Torvik setup}. Note that
the displacement $y\left( t\right) $ of the plate is referred to the
equilibrium point $\ell _{0}$ at which the weight of the plate is
compensated by the force exerted by the spring and the buoyancy force
experimented by the plate.

The motion of the plate is governed by the following fractional-order
differential equation \cite[Eqn. 8.18]{Podlubny}:%
\begin{equation}
a\,\,y^{\prime \prime }\left( t\right) +b\,\,_{0}^{\ast }D_{t}^{3/2}y\left(
t\right) +c\,y\left( t\right) =f\left( t\right) ,\quad t>0,
\label{BT_equation}
\end{equation}%
where $a=m$, $b=2A\sqrt{\mu \rho }$, $c=K$ and $_{0}^{\ast
}D_{t}^{3/2}y\left( t\right) $ denotes the Caputo fractional derivative of
order $3/2$. Next, we recall the definition of the Caputo fractional
derivative.

\begin{definition}
For $m=1,2,\ldots$, the Caputo fractional derivative of order $\alpha >0$ is defined as \cite[%
Eqn. 1.17b]{MainardiBook}:%
\begin{equation}
_{0}^{\ast }D_{t}^{\alpha }g\left( t\right) =\left\{
\begin{array}{ll}
\displaystyle%
\frac{1}{\Gamma \left( m-\alpha \right) }\int_{0}^{t}\frac{g^{\left(
m\right) }\left( \tau \right) }{\left( t-\tau \right) ^{\alpha +1-m}}d\tau ,
& m-1<\alpha <m, \\
\displaystyle \frac{d^{m}}{dt^{m}}g\left( t\right) , & \alpha =m.%
\end{array}%
\right.  \label{Caputo_def}
\end{equation}
\end{definition}

It is worth noting that the plate-fluid system must be in an equilibrium
state in order to derive (\ref{BT_equation}), thus the initial velocity of the
plate must be zero, i.e. $y^{\prime }\left( 0\right) =0$. In the original
work of Torvik and Bagley \cite{BangleyTorvik1984}, they considered
homogeneous initial conditions:
\begin{equation}
y\left( 0\right) =y^{\prime }\left( 0\right) =0.  \label{Initial_homo}
\end{equation}

Here, we generalize the Bagley-Torvik equation (\ref{BT_equation})\ for
non-homogeneous initial conditions $y\left( 0\right) $, $y^{\prime }\left(
0\right) $. In this case, we can write the solution as
\begin{equation}
y\left( t\right) =y_{f}\left( t\right) +y_{c}\left( t\right) ,
\label{Intro_1}
\end{equation}%
where $y_{f}\left( t\right) $ takes into account the effect of the external
force $f\left( t\right) $ exerted to the plate, and $y_{c}\left( t\right) $
takes into account the effect of non-homogeneous initial conditions. Note
that if $f\left( t\right) =0$ and $y\left( 0\right) =y^{\prime }\left(
0\right) =0$, the solution of (\ref{BT_equation})\ is $y\left( t\right) =0$,
for $t>0$. Physically speaking, this means that the plate will remain at the
equilibrium point, i.e. the initial position.

Since the Bagley-Torvik equation plays a vital role in many problems of
applied science and engineering, a large number of researchers show interest
in its solution. In Section \ref{Section: Preliminaries},\ we explain the
most common analytical and numerical approaches to solve the Bagley-Torvik
equation and its generalizations. Also in Section \ref{Section: Preliminaries}, we present some properties of the Mittag-Leffler function
that will be needed throughout the paper. The reader will find many more references
in \cite{Zafar} regarding the great variety of approaches found in the
literature in order to solve the Bagley-Torvik equation.

It is worth noting that a fractional-order differential equation similar to (%
\ref{BT_equation}) arises in fluid mechanics when we study the motion of a
rigid sphere in the gravity field, starting from rest in a quiescent fluid.
This problem was first treated by Basset \cite{Basset} and since then is
called Basset problem. The generalization of the Basset problem and its
analytical solution via Laplace transform method is given in \cite%
{MainardiBasset}. We will follow a similar approach in order to analytically
solve (\ref{BT_equation}) in this paper.

The goal of this paper is to solve the Bagley-Torvik equation given in (\ref%
{BT_equation})\ for given constants $a,b,c\in
\mathbb{R}
$ and arbitrary external force $f\left( t\right) $ with initial conditions $%
y\left( 0\right) ,y^{\prime }\left( 0\right) \in
\mathbb{R}
$ in closed-form. For this purpose, we derive a new Laplace transform in
Section \ref{Section: Laplace inversion}. In Section \ref{Section: BT
general solution}, we apply this new Laplace transform to obtain the general
solution of (\ref{BT_equation}), where $y_{c}\left( t\right) $ is given as a
linear combination of exponentials with error functions;\ and where $%
y_{f}\left( t\right) $ is given in terms of a convolution integral. In the
next two sections, we calculate $y_{f}\left( t\right) $ for particular forms
of $f\left( t\right) $. Thereby, in Section \ref{Section: potential solution}%
, we consider $f\left( t\right) $ as a linear combination of potential
functions; and in Section \ref{Section: sinusoidal solution}, as a
sinusoidal function. In addition we present the asymptotic behaviour of the
solution for $t\rightarrow 0^{+}$ and $t\rightarrow \infty$ in Section \ref{Section: asymptotic}.
In Section \ref{Section: Numerical examples}, we
present some numerical examples in order to compare the behaviour of our
solution to other approaches found in the literature. We collect our
conclusions in Section \ref{Section: Conclusions}. Finally, the calculation
of some auxiliary integrals 
is included in Appendix \ref{secA1}.

\section{Preliminaries \label{Section: Preliminaries}}
As aforementioned in the Introduction, this section is devoted to presenting some properties of the Mittag-Leffler function
that will be used throughout the paper, as well as to explain the
most common analytical and numerical approaches in existing literature in order to solve the Bagley-Torvik
equation.
\subsection{Some properties of the Mittag-Leffler function}

\begin{definition}
The one-parameter Mittag-Leffler function is defined as \cite[Eqn. 3.1.1]{GorenfloML}:
\begin{equation}
\mathrm{E}_{\alpha }\left( z\right) =\sum_{k=0}^{\infty }\frac{z^{k}}{\Gamma
\left( \alpha k+1\right) },\quad \alpha \in \mathbb{C},  \label{ML_1parameter_def}
\end{equation}
where $\Gamma \left( z\right) $ denotes the gamma function \cite[Eqn. 43:3:1]{Atlas}.
\end{definition}

\begin{definition}
The two parameter Mittag-Leffler function is defined as \cite[Eqn. 4.1.1]{GorenfloML}:
\begin{equation}
\mathrm{E}_{\alpha ,\beta }\left( z\right) =\sum_{k=0}^{\infty }\frac{z^{k}}{%
\Gamma \left( \alpha k+\beta \right) },\quad \mathrm{Re}\,\alpha >0 ,\beta \in \mathbb{C}.
\label{ML_2parameter_def}
\end{equation}
\end{definition}

\begin{remark}
Note that the one-parameter Mittag-Leffler function is a special case of the
two-parameter Mittag-Leffler function, since $\mathrm{E}_{\alpha }\left(
z\right) =\mathrm{E}_{\alpha ,1}\left( z\right) $.
\end{remark}

In existing literature \cite[Eqn. 3.2.4]{GorenfloML}, 
we found the reduction formula%
\begin{equation}
\mathrm{E}_{1/2,1}\left( z\right) =\mathrm{E}_{1/2}\left( z\right) =\exp
\left( z^{2}\right) \,\mathrm{erfc}\left( -z\right) .  \label{ML(1/2,1)}
\end{equation}%
Also, we found the property \cite[Eqn. 4.2.3]{GorenfloML}:%
\begin{equation}
\mathrm{E}_{\alpha ,\beta }\left( z\right) =\frac{1}{\Gamma \left( \beta
\right) }+z\,\mathrm{E}_{\alpha ,\alpha +\beta }\left( z\right) .
\label{Recursive_ML}
\end{equation}%
From (\ref{ML(1/2,1)})\ and (\ref{Recursive_ML}), and knowing that $\Gamma
\left( \frac{1}{2}\right) =\sqrt{\pi }$ and $\Gamma \left( 1\right) =1$ \cite%
[Sect. 43:7]{Atlas}, we obtain%
\begin{equation}
\mathrm{E}_{1/2,1/2}\left( z\right) =\frac{1}{\sqrt{\pi }}+z\exp \left(
z^{2}\right) \,\mathrm{erfc}\left( -z\right) ,  \label{ML(1/2,1/2)}
\end{equation}%
and%
\begin{equation}
\mathrm{E}_{1/2,3/2}\left( z\right) =\frac{\exp \left( z^{2}\right) \,%
\mathrm{erfc}\left( -z\right) -1}{z}.  \label{ML(1/2,3/2)}
\end{equation}%
From (\ref{ML(1/2,1/2)})\ and (\ref{Recursive_ML}), and knowing that $%
1/\Gamma \left( 0\right) =0$ \cite[Sect. 43:7]{Atlas}, we obtain%
\begin{equation}
\mathrm{E}_{1/2,0}\left( z\right) =z\left( \frac{1}{\sqrt{\pi }}+z\exp
\left( z^{2}\right) \,\mathrm{erfc}\left( -z\right) \right) .
\label{ML(1/2,0)}
\end{equation}%
Finally, from (\ref{ML(1/2,0)})\ and (\ref{Recursive_ML}), and knowing that $%
\Gamma \left( -\frac{1}{2}\right) =-2\sqrt{\pi }$ \cite[Sect. 43:7]{Atlas}, we obtain%
\begin{equation}
\mathrm{E}_{1/2,-1/2}\left( z\right) =\frac{1}{\sqrt{\pi }}\left( z^{2}-%
\frac{1}{2}\right) +z^{3}\exp \left( z^{2}\right) \,\mathrm{erfc}\left(
-z\right) .  \label{ML(1/2,-1/2)}
\end{equation}

\subsection{Analytical solutions}

The analytical solution of (\ref{BT_equation})\ given in the literature for
homogeneous initial conditions (\ref{Initial_homo}) reads as \cite[Eqns.
8.26-27]{Podlubny}, :%
\begin{equation}
y\left( t\right) =y_{f}\left( t\right) =\int_{0}^{t}f\left( \tau \right)
G\left( t-\tau \right) d\tau ,  \label{yf(t)_Podlubny}
\end{equation}%
being%
\begin{equation}
G\left( t\right) =\frac{1}{a}\sum_{k=0}^{\infty }\frac{\left( -1\right) ^{k}%
}{k!}\left( \frac{c}{a}\right) ^{k}t^{2k+1}\,\mathrm{E}_{1/2,2+3k/2}^{\left(
k\right) }\left( -\frac{b}{a}\sqrt{t}\right) ,  \label{G(t)_Podlubny}
\end{equation}%
and $\mathrm{E}_{\lambda ,\mu }^{\left( k\right) }\left( z\right) $ the $%
k $-th derivative of the two-parameter Mittag-Leffler function (\ref{ML_2parameter_def}). 
\begin{equation}
\mathrm{E}_{\lambda ,\mu }^{\left( k\right) }\left( z\right) =\frac{d^{k}}{%
dz^{k}}\mathrm{E}_{\lambda ,\mu }\left( z\right) =\sum_{j=0}^{\infty }\frac{%
\left( j+k\right) !\,z^{j}}{j!\,\Gamma \left( \lambda j+\lambda k+\mu
\right) }.  \label{Dk_ML_def}
\end{equation}%

By using the Adomian decomposition, the solution given in (\ref%
{yf(t)_Podlubny})-(\ref{Dk_ML_def}) is obtained in \cite{Adomian}\ for the
particular case $f\left( t\right) =8\left[ H\left( t\right) - H\left( t-1\right) \right] $ where $ H\left( t\right) $ denotes the
Heaviside function. By using the Laplace transform, the solution given in (%
\ref{yf(t)_Podlubny})-(\ref{Dk_ML_def}) is obtained in \cite[Example 8.2]%
{Diethelm}\ for the particular case $f\left( t\right) =\sin t$.

In \cite{Pang}, we found the solution of (\ref{BT_equation})\ for
non-homogeneous initial conditions $y\left( 0\right) $, $y^{\prime }\left(
0\right) $, where $y_{f}\left( t\right) $ is given by (\ref{yf(t)_Podlubny}%
)-(\ref{Dk_ML_def})\ and
\begin{eqnarray}
y_{c}\left( t\right) &=&y\left( 0\right) \sum_{k=0}^{\infty }\frac{\left(
-1\right) ^{k}}{k!}\left( \frac{c}{a}\right) ^{k}t^{2k}\,\mathrm{E}%
_{1/2,1+3k/2}^{\left( k\right) }\left( -\frac{b}{a}\sqrt{t}\right)
\label{yc_Pang} \\
&&+y\left( 0\right) \frac{b}{a}\sum_{k=0}^{\infty }\frac{\left( -1\right)
^{k}}{k!}\left( \frac{c}{a}\right) ^{k}t^{2k+1/2}\,\mathrm{E}_{1/2,3+3\left(
k-1\right) /2}^{\left( k\right) }\left( -\frac{b}{a}\sqrt{t}\right)  \notag
\\
&&+y^{\prime }\left( 0\right) \frac{b}{a}\sum_{k=0}^{\infty }\frac{\left(
-1\right) ^{k}}{k!}\left( \frac{c}{a}\right) ^{k}t^{2k+3/2}\,\mathrm{E}%
_{1/2,4+3\left( k-1\right) /2}^{\left( k\right) }\left( -\frac{b}{a}\sqrt{t}%
\right)  \notag \\
&&+y^{\prime }\left( 0\right) \sum_{k=0}^{\infty }\frac{\left( -1\right) ^{k}%
}{k!}\left( \frac{c}{a}\right) ^{k}t^{2k+1}\,\mathrm{E}_{1/2,2+3k/2}^{\left(
k\right) }\left( -\frac{b}{a}\sqrt{t}\right) .  \notag
\end{eqnarray}

Other approach to obtain an analytical solution for non-homogeneous initial
conditions is given in \cite{Arora}. However, this method is restricted to
functions of the following type:%
\begin{equation}
f\left( t\right) =\sum_{k=0}^{n}a_{k}\,t^{k/2}.  \label{f(t)_Arora}
\end{equation}%
The solution is given in series form as\
\begin{equation}
y\left( t\right) =\sum_{k=0}^{\infty }d_{k}\,t^{k/2},  \label{y(t)_Arora}
\end{equation}%
where $d_{0}=y\left( 0\right) $, $d_{1}=0$, $d_{2}=y^{\prime }\left(
0\right) $, $d_{3}=0$, and%
\begin{eqnarray}
d_{k+4} &=&\frac{\Gamma \left( \frac{k}{2}+1\right) \left[ \phi \left(
k\right) -c\,d_{k}\right] -b\,\Gamma \left( \frac{k+5}{2}\right) \,d_{k+3}}{%
a\,\Gamma \left( \frac{k}{2}+3\right) },  \label{Recursive_Arora} \\
k &=&0,1,2,\ldots ,  \notag
\end{eqnarray}%
with%
\begin{equation}
\phi \left( k\right) =\sum_{m=0}^{n}a_{m}\,\delta _{k-m,0},
\label{Phi(k)_Arora}
\end{equation}%
being $\delta _{n,m}$ the Kronecker delta function \cite[Eqn. 9:13:2]{Atlas}.

Other generalizations of the Bagley-Torvik equation with non-homogeneous
initial conditions can be found in \cite{Zafar}, i.e.
\begin{equation}
a\,\,_{0}^{\ast }D_{t}^{\beta }y\left( t\right) +b\,\,_{0}^{\ast
}D_{t}^{\alpha +1}y\left( t\right) +c\,y\left( t\right) =f\left( t\right)
,\quad t>0,  \label{BTE_alfa_beta}
\end{equation}%
where $1<\beta <2$ and $0<\alpha <1$. The analytical solution of (\ref%
{BTE_alfa_beta})\ is given in terms of the Lorenzo-Hartley \textquotedblleft
generalized $G$ function\textquotedblright\ \cite[Eqn. 101]{LorenzoHartley}:%
\begin{equation}
G_{a,b,c}\left( z,t\right) =\sum_{j=0}^{\infty }\frac{z^{j}\Gamma \left(
c+j\right) \,t^{\left( c+j\right) a-b-1}}{\Gamma \left( c\right) \,\Gamma
\left( j+1\right) \,\Gamma \left( \left( c+j\right) a-b\right) }.
\label{G_Hartley_def}
\end{equation}

\subsection{Numerical solutions}

Next, we describe the numerical method given by Podlubny in \cite[Sect. 8.3.2]{Podlubny}
to solve the Bagley-Torvik equation (\ref{BT_equation}). First, we
discretize the solution $y\left( t\right) $ as well as the function $f\left(
t\right) $ with a time step $h$ as
\begin{eqnarray}
y\left( t\right) &\approx &y_{m}=y\left( mh\right) ,  \label{ym_def} \\
f\left( t\right) &\approx &f_{m}=f\left( mh\right) ,  \label{fm_def}
\end{eqnarray}%
where $m=0,1,2,\ldots $ and $t$ is the point at which the function is
evaluated, being $m=t/h$ the number of nodes. According to \cite[Eqns.
7.9\&7.22]{Podlubny}, the first-order approximation for the $\alpha $-th
derivative is
\begin{equation}
_{0}D_{t}^{\alpha }\,g\left( t\right) \approx h^{-\alpha
}\sum_{j=0}^{m}w_{j}^{\left( \alpha \right) }g\left( t-jh\right) ,
\label{0_D_t_(a)_approx}
\end{equation}%
where
\begin{equation}
w_{j}^{\left( \alpha \right) }=\left( -1\right) ^{j}\binom{\alpha }{j}.
\label{wj_def}
\end{equation}%
Consequently, from (\ref{ym_def})-(\ref{wj_def}), and taking into account
that $w_{0}^{\left( \alpha \right) }=1$, the first order approximation of
the Bagley-Torvik equation (\ref{BT_equation}), is
\begin{equation}
a\,\underset{\approx \,y^{\prime \prime }\left( t\right) }{\underbrace{\frac{%
y_{m}-2\,y_{m-1}+y_{m-2}}{h^{2}}}}+b\,\underset{\approx
\,_{0}D_{t}^{3/2}y\left( t\right) }{\underbrace{h^{-3/2}%
\sum_{j=0}^{m}w_{j}^{\left( \alpha \right) }y_{m-j}}}+c\,y_{m}=f_{m},
\label{Podlubny_1}
\end{equation}%
thus (it is worth noting that there is an error in Podlubny's treatise \cite[Eqn. 8.25]{Podlubny}):%
\begin{equation}
y_{m}=\frac{h^{2}f_{m}+a\left( 2\,y_{m-1}-y_{m-2}\right) -b\sqrt{h}%
\sum_{j=1}^{m}w_{j}^{\left( \alpha \right) }y_{m-j}}{a+b\sqrt{h}+c\,h^{2}}.
\label{Podlubny_2}
\end{equation}%
where $m=2,3,\ldots $ and the initial conditions are approximated by
\begin{equation}
y_{0}=y\left( 0\right) ,\quad y_{1}=y\left( 0\right) +h\,y^{\prime }\left(
0\right) .  \label{Podlubny_initial}
\end{equation}

Diethelm reformulate (\ref{BT_equation}) into the following system of
fractional differential equations:%
\begin{equation}
\begin{array}{l}
_{0}^{\ast }D_{t}^{1/2}y_{1}\left( t\right) =y_{2}\left( t\right) , \\
_{0}^{\ast }D_{t}^{1/2}y_{2}\left( t\right) =y_{3}\left( t\right) , \\
_{0}^{\ast }D_{t}^{1/2}y_{3}\left( t\right) =y_{4}\left( t\right) , \\
\displaystyle%
_{0}^{\ast }D_{t}^{1/2}y_{4}\left( t\right) =\frac{1}{a}\left[
-c\,y_{1}\left( t\right) -b\,y_{4}\left( t\right) +f\left( t\right) \right] ,%
\end{array}
\label{Diethelm_1}
\end{equation}%
with $y_{1}\left( 0\right) =y\left( 0\right) $, $y_{2}\left( 0\right) =0$, $%
y_{3}\left( 0\right) =y^{\prime }\left( 0\right) $, $y_{4}\left( 0\right) =0$%
, and then numerically solve the model with the Adam predictor-corrector
approach \cite{DiethelmNumerical}. Very recently, in \cite{Bonis}, the
authors used a new approximation to compute the Caputo fractional derivative
in order to numerically solve (\ref{BT_equation}) with non-homogeneous
initial conditions. Another generalization of the Bagley-Torvik equation
with non-homogeneous initial conditions is given in \cite{YubasiNumerical}\
\begin{equation}
a\left( t\right) \,\,y^{\prime \prime }\left( t\right) +b\left( t\right)
\,_{0}^{\ast }D_{t}^{3/2}y\left( t\right) +c\left( t\right) \,y\left(
t\right) =f\left( t\right) ,  \label{Yubasi_1}
\end{equation}%
which is numerically solved by using a Bessel collocation method. A further
generalization of the Bagley-Torvik equation is found in \cite%
{GulsuNumerical}, where the authors numerically solve the following
fractional differential equation:%
\begin{equation}
\begin{array}{l}
\displaystyle%
\sum_{k=0}^{m}a_{k}\left( t\right) \,_{0}^{\ast }D_{t}^{k\alpha }y\left(
t\right) \,=f\left( t\right) , \\
a\leq t\leq b,\ n-1\leq m\alpha <n,%
\end{array}
\label{Gulsu_1}
\end{equation}%
\ with initial conditions%
\begin{equation*}
_{0}^{\ast }D_{t}^{i}y\left( c\right) ,\quad i=0,1,\ldots ,n-1,\ a\leq c\leq
b.
\end{equation*}

\section{A useful inverse Laplace transform \label{Section: Laplace
inversion}}

%

Let us denote the Laplace transform of a function $g\left(t \right)$ for $t>0$ as $G\left( s\right) =\mathcal{L}\left[ g\left( t\right) ;s\right]$, thus the inverse Laplace transform is denoted as $\mathcal{L}^{-1}\left[ G\left( s\right) ;t\right] =g\left( t\right)$. We know that the Laplace transform of the Caputo fractional derivative is \cite[Eqn. 1.27]%
{MainardiBook}
\begin{equation}
\mathcal{L}\left[ _{0}^{\ast }D_{t}^{\alpha }g\left( t\right) ;s\right]
=s^{\alpha }\mathcal{L}\left[ g\left( t\right) ;s\right] -\sum_{k=0}^{m}s^{%
\alpha -1-k}g\left( 0\right) .  \label{Laplace_Caputo}
\end{equation}
For $\alpha =m$, (\ref{Laplace_Caputo})\ is reduced to (see also \cite[%
Theorem 2.12]{Schiff})%
\begin{equation}
\mathcal{L}\left[ g^{\left( m\right) }\left( t\right) ;s\right] =s^{m}%
\mathcal{L}\left[ g\left( t\right) ;s\right] -\sum_{k=0}^{m}s^{m-1-k}g\left(
0\right) .  \label{Laplace_derivative_m}
\end{equation}

Now, we want to generalize the recent inverse Laplace transform calculated
in \cite{NovelLaplaceJL}. For this purpose, let us calculate the inverse
Laplace transform of the function:%
\begin{equation}
G_{m,p/q,\lambda }\left( s\right) =\frac{s^{\lambda }}{a\,s^{m}+b\,s^{p/q}+c}%
,  \label{Laplace_1}
\end{equation}%
where $p/q\in
\mathbb{Q}
^{+}$, $m\in
\mathbb{N}
$, $\lambda \in
\mathbb{R}
$. Thereby,
\begin{eqnarray}
g_{m,p/q,\lambda }\left( t\right) &=&\mathcal{L}^{-1}\left[ G_{m,p/q,\lambda
}\left( s\right) ;t\right]  \notag \\
&=&\mathcal{L}^{-1}\left[ \frac{s^{\lambda }}{a\,s^{m}+b\,s^{p/q}+c};t\right]
.  \label{Laplace_2}
\end{eqnarray}%
Perform the change of variables $r=s^{1/q}$, to obtain%
\begin{equation}
g_{m,p/q,\lambda }\left( t\right) =\mathcal{L}^{-1}\left[ \frac{r^{q\lambda }%
}{P\left( r\right) };t\right] ,  \label{g_m(t)_a}
\end{equation}%
where we define the polynomial
\begin{equation}
P\left( r\right) =a\,r^{mq}+b\,r^{p}+c.  \label{P(r)_general_def}
\end{equation}%
If $P\left( r\right) $ has $N=\max \left( mq,p\right) $ different roots $%
r_{k}$ $\left( k=1,\ldots ,N\right) $, according to \cite[Eqn. 17:13:10]%
{Atlas}, we can rewrite (\ref{g_m(t)_a})\ as
\begin{eqnarray}
g_{m,p/q,\lambda }\left( t\right) &=&\mathcal{L}^{-1}\left[ \sum_{k=1}^{N}%
\frac{r^{q\lambda }}{P^{\prime }\left( r_{k}\right) \left( r-r_{k}\right) };t%
\right]  \notag \\
&=&\sum_{k=1}^{N}\frac{1}{P^{\prime }\left( r_{k}\right) }\mathcal{L}^{-1}%
\left[ \frac{s^{\lambda }}{s^{1/q}-r_{k}};t\right] .  \label{Laplace_3}
\end{eqnarray}%
Applying the inverse Laplace transform \cite[Eqn. 45:14:4]{Atlas}:%
\begin{equation}
\mathcal{L}^{-1}\left[ \frac{s^{\mu -\nu }}{s^{\mu }-\alpha };t\right]
=t^{\mu -1}\mathrm{E}_{\mu ,\nu }\left( \alpha \,t^{\mu }\right) ,
\label{Laplace_4}
\end{equation}%
where $\mathrm{E}_{\mu ,\nu }\left( z\right) $ denotes the two-parameter
Mittag-Leffler function, we arrive at%
\begin{equation}
g_{m,p/q,\lambda }\left( t\right) =t^{1/q-\lambda -1}\sum_{k=1}^{N}\frac{%
\mathrm{E}_{1/q,1/q-\lambda }\left( r_{k}\,t^{1/q}\right) }{P^{\prime
}\left( r_{k}\right) }.  \label{Laplace_5}
\end{equation}

We summarize the above calculation as follows.

\begin{theorem}
For $p/q\in
\mathbb{Q}
^{+}$, $m\in
\mathbb{N}
$, $\lambda \in
\mathbb{R}
$, the following inverse Laplace transform holds true:\
\begin{equation}
\mathcal{L}^{-1}\left[ \frac{s^{\lambda }}{a\,s^{m}+b\,s^{p/q}+c};t\right]
=t^{1/q-\lambda -1}\sum_{k=1}^{N}\frac{\mathrm{E}_{1/q,1/q-\lambda }\left(
r_{k}\,t^{1/q}\right) }{P^{\prime }\left( r_{k}\right) },
\label{Inverse_Laplace_transform}
\end{equation}%
where $r_{k}$ $\left( k=1,\ldots ,N\right) $ are the $N=\max \left(
mq,p\right) $ different roots of the polynomial:
\begin{equation}
P\left( r\right) =a\,r^{mq}+b\,r^{p}+c.  \label{Poly_def}
\end{equation}
\end{theorem}

\section{The general solution of the Bagley-Torvik equation \label{Section:
BT general solution}}

Apply the Laplace transform to the Bagley-Torvik equation (\ref{BT_equation}%
), taking into account the properties of the Laplace transform given in (\ref%
{Laplace_Caputo}) and (\ref{Laplace_derivative_m}), to obtain
\begin{eqnarray}
&&a\,\left[ s^{2}Y\left( s\right) -s\,y\left( 0\right) -y^{\prime }\left(
0\right) \right]  \label{BT_1} \\
+ &&b\,\left[ s^{3/2}Y\left( s\right) -s^{1/2}y\left( 0\right)
-s^{-1/2}y^{\prime }\left( 0\right) \right] +c\,Y\left( s\right) =F\left(
s\right) ,  \notag
\end{eqnarray}%
thus%
\begin{eqnarray}
Y\left( s\right) &=&\frac{F\left( s\right) }{a\,s^{2}+b\,s^{3/2}+c}
\label{Y(s)_sol} \\
&&+y\left( 0\right) \left[ \frac{a\,s}{a\,s^{2}+b\,s^{3/2}+c}+\frac{%
b\,s^{1/2}}{a\,s^{2}+b\,s^{3/2}+c}\right]  \notag \\
&&+y^{\prime }\left( 0\right) \left[ \frac{a}{a\,s^{2}+b\,s^{3/2}+c}+\frac{%
b\,s^{-1/2}}{a\,s^{2}+b\,s^{3/2}+c}\right] .  \notag
\end{eqnarray}%
Now, perform the inverse Laplace transform in order to write the solution as%
\begin{equation}
y\left( t\right) =y_{f}\left( t\right) +y_{c}\left( t\right) ,
\label{y=yf+yc_def}
\end{equation}%
where we define%
\begin{equation}
y_{f}\left( t\right) =\mathcal{L}^{-1}\left[ \frac{F\left( s\right) }{%
a\,s^{2}+b\,s^{3/2}+c};t\right] ,  \label{yf_def}
\end{equation}%
and%
\begin{eqnarray}
&&y_{c}\left( t\right)  \label{yc_def} \\
&=&y\left( 0\right) \left\{ a\,\mathcal{L}^{-1}\left[ \frac{s}{%
a\,s^{2}+b\,s^{3/2}+c};t\right] +b\,\mathcal{L}^{-1}\left[ \frac{s^{1/2}}{%
a\,s^{2}+b\,s^{3/2}+c};t\right] \right\}  \notag \\
&&+y^{\prime }\left( 0\right) \left\{ a\,\mathcal{L}^{-1}\left[ \frac{1}{%
a\,s^{2}+b\,s^{3/2}+c};t\right] +b\,\mathcal{L}^{-1}\left[ \frac{s^{-1/2}}{%
a\,s^{2}+b\,s^{3/2}+c};t\right] \right\} .  \notag
\end{eqnarray}

Let us calculate the inverse Laplace transforms given in (\ref{yc_def}).
First, apply (\ref{Inverse_Laplace_transform}) for $\lambda =-\frac{1}{2},0,%
\frac{1}{2},1$, and $p=3$, $q=2$, $m=2$, to obtain%
\begin{eqnarray}
\mathcal{L}^{-1}\left[ \frac{s^{-1/2}}{a\,s^{2}+b\,s^{3/2}+c};t\right]
&=&\sum_{k=1}^{4}\frac{\mathrm{E}_{1/2,1}\left( r_{k}\,\sqrt{t}\right) }{%
\left( 4a\,r_{k}+3b\right) r_{k}^{2}},  \label{Inverse_-1/2} \\
\mathcal{L}^{-1}\left[ \frac{1}{a\,s^{2}+b\,s^{3/2}+c};t\right] &=&\frac{1}{%
\sqrt{t}}\sum_{k=1}^{4}\frac{\mathrm{E}_{1/2,1/2}\left( r_{k}\,\sqrt{t}%
\right) }{\left( 4a\,r_{k}+3b\right) r_{k}^{2}},  \label{Inverse_0} \\
\mathcal{L}^{-1}\left[ \frac{s^{1/2}}{a\,s^{2}+b\,s^{3/2}+c};t\right] &=&%
\frac{1}{t}\sum_{k=1}^{4}\frac{\mathrm{E}_{1/2,0}\left( r_{k}\,\sqrt{t}%
\right) }{\left( 4a\,r_{k}+3b\right) r_{k}^{2}},  \label{Inverse_1/2} \\
\mathcal{L}^{-1}\left[ \frac{s}{a\,s^{2}+b\,s^{3/2}+c};t\right] &=&\frac{1}{t%
\sqrt{t}}\sum_{k=1}^{4}\frac{\mathrm{E}_{1/2,-1/2}\left( r_{k}\,\sqrt{t}%
\right) }{\left( 4a\,r_{k}+3b\right) r_{k}^{2}},  \label{Inverse_1}
\end{eqnarray}%
where $r_{k}$ $\left( k=1,2,3,4\right) $ are the $4$ different roots of the
polynomial:
\begin{equation}
P\left( r\right) =a\,r^{4}+b\,r^{3}+c.  \label{quartic_polynomial}
\end{equation}%
In order to obtain the roots of (\ref{quartic_polynomial}), define:%
\begin{equation}
\begin{array}{l}
\displaystyle%
\beta =\frac{b}{a},\ \gamma =\frac{c}{a}, \\
\displaystyle%
\delta =\left( 9\beta ^{2}+\sqrt{81\beta ^{4}-768\gamma }\right) \gamma , \\
\displaystyle%
R=\frac{1}{2}\sqrt{\left( \frac{\beta }{2}\right) ^{2}+\left( \frac{\delta }{%
18}\right) ^{1/3}+\frac{8\gamma }{\left( 12\delta \right) ^{1/3}}}, \\
\displaystyle%
T_{\pm }=\frac{1}{2}\sqrt{\frac{\beta ^{2}}{2}-4\gamma \left( \frac{2}{%
3\delta }\right) ^{1/3}-\left( \frac{\delta }{18}\right) ^{1/3}\mp \frac{%
\beta ^{3}}{8R}}.%
\end{array}
\label{Parameters_quartic}
\end{equation}%
\ thus, according to \cite[Sect. 9]{Spiegel}, we have%
\begin{equation}
\begin{array}{c}
\displaystyle%
r_{1,2}=-\frac{\beta }{4}\pm R+T_{\pm }, \\
\displaystyle%
r_{3,4}=-\frac{\beta }{4}\pm R-T_{\pm }.%
\end{array}
\label{rk_roots}
\end{equation}
According to the reduction formulas given in (\ref{ML(1/2,1)})-(\ref{ML(1/2,-1/2)}), 
we rewrite (\ref{Inverse_-1/2})-(\ref{Inverse_1}) as
\begin{eqnarray}
\mathcal{L}^{-1}\left[ \frac{s^{-1/2}}{a\,s^{2}+b\,s^{3/2}+c};t\right]
&=&\sum_{k=1}^{4}\frac{\exp \left( r_{k}^{2}\,t\right) \,\mathrm{erfc}\left(
-r_{k}\,\sqrt{t}\right) }{\left( 4a\,r_{k}+3b\right) r_{k}^{2}},
\label{Inv_-1/2} \\
\mathcal{L}^{-1}\left[ \frac{1}{a\,s^{2}+b\,s^{3/2}+c};t\right] &=&\frac{%
A_{-2}}{\sqrt{\pi t}}+\sum_{k=1}^{4}\frac{\exp \left( r_{k}^{2}\,t\right) \,%
\mathrm{erfc}\left( -r_{k}\,\sqrt{t}\right) }{\left( 4a\,r_{k}+3b\right)
r_{k}},  \label{Inv_0} \\
\mathcal{L}^{-1}\left[ \frac{s^{1/2}}{a\,s^{2}+b\,s^{3/2}+c};t\right] &=&%
\frac{A_{-1}}{\sqrt{\pi t}}+\sum_{k=1}^{4}\frac{\exp \left(
r_{k}^{2}\,t\right) \,\mathrm{erfc}\left( -r_{k}\,\sqrt{t}\right) }{%
4a\,r_{k}+3b},  \label{Inv_1/2} \\
\mathcal{L}^{-1}\left[ \frac{s}{a\,s^{2}+b\,s^{3/2}+c};t\right] &=&\frac{%
A_{0}}{\sqrt{\pi t}}-\frac{A_{-2}}{2t\sqrt{\pi t}}  \label{Inv_1} \\
&&+\sum_{k=1}^{4}\frac{r_{k}\exp \left( r_{k}^{2}\,t\right) \,\mathrm{erfc}%
\left( -r_{k}\,\sqrt{t}\right) }{4a\,r_{k}+3b},  \notag
\end{eqnarray}%
where we have defined%
\begin{equation}
A_{\ell }=\sum_{k=1}^{4}\frac{r_{k}^{\ell }}{4a\,r_{k}+3b}.  \label{A_l_def}
\end{equation}
Note that the inverse Laplace transforms given in (\ref{Inv_-1/2})-(\ref%
{Inv_1})\ should be non-divergent for $t=0$. Consequently,
\begin{equation}
A_{0}=A_{-1}=A_{-2}=0,  \label{A_l=0}
\end{equation}%
as we can numerically check from (\ref{A_l_def}). Therefore, (\ref{Inv_-1/2}%
)-(\ref{Inv_1}) are reduced to%
\begin{eqnarray}
\mathcal{L}^{-1}\left[ \frac{s^{-1/2}}{a\,s^{2}+b\,s^{3/2}+c};t\right]
&=&\sum_{k=1}^{4}\frac{\exp \left( r_{k}^{2}\,t\right) \,\mathrm{erfc}\left(
-r_{k}\,\sqrt{t}\right) }{\left( 4a\,r_{k}+3b\right) r_{k}^{2}},
\label{Inverse_s^(-1/2)} \\
\mathcal{L}^{-1}\left[ \frac{1}{a\,s^{2}+b\,s^{3/2}+c};t\right]
&=&\sum_{k=1}^{4}\frac{\exp \left( r_{k}^{2}\,t\right) \,\mathrm{erfc}\left(
-r_{k}\,\sqrt{t}\right) }{\left( 4a\,r_{k}+3b\right) r_{k}},
\label{Inverse_s^0} \\
\mathcal{L}^{-1}\left[ \frac{s^{1/2}}{a\,s^{2}+b\,s^{3/2}+c};t\right]
&=&\sum_{k=1}^{4}\frac{\exp \left( r_{k}^{2}\,t\right) \,\mathrm{erfc}\left(
-r_{k}\,\sqrt{t}\right) }{4a\,r_{k}+3b},  \label{Inverse_s^(1/2)} \\
\mathcal{L}^{-1}\left[ \frac{s}{a\,s^{2}+b\,s^{3/2}+c};t\right]
&=&\sum_{k=1}^{4}\frac{r_{k}\exp \left( r_{k}^{2}\,t\right) \,\mathrm{erfc}%
\left( -r_{k}\,\sqrt{t}\right) }{4a\,r_{k}+3b}.  \label{Inverse_s^1}
\end{eqnarray}%
Insert (\ref{Inverse_s^(-1/2)})-(\ref{Inverse_s^1})\ into (\ref{yc_def}), to
arrive at%
\begin{equation}
y_{c}\left( t\right) =\sum_{k=1}^{4}\frac{\left( a\,r_{k}+b\right) \left[
r_{k}^{2}\,y\left( 0\right) +y^{\prime }\left( 0\right) \right] }{\left(
4a\,r_{k}+3b\right) r_{k}^{2}}\exp \left( r_{k}^{2}\,t\right) \,\mathrm{erfc}%
\left( -r_{k}\,\sqrt{t}\right) .  \label{yc(t)_resultado}
\end{equation}

Second, let us calculate the inverse Laplace transform given in (\ref{yf_def}%
), taking into account the result obtained in (\ref{Inverse_s^0}),
\begin{eqnarray}
y_{f}\left( t\right) &=&\mathcal{L}^{-1}\left[ F\left( s\right) \ \frac{1}{%
a\,s^{2}+b\,s^{3/2}+c};t\right]  \notag \\
&=&\mathcal{L}^{-1}\left\{ \mathcal{L}\left[ f\left( t\right) ;s\right] \
\mathcal{L}\left[ \mathcal{L}^{-1}\left( \frac{1}{a\,s^{2}+b\,s^{3/2}+c}%
;t\right) ;s\right] \right\}  \notag \\
&=&\mathcal{L}^{-1}\left\{ \mathcal{L}\left[ f\left( t\right) ;s\right] \
\mathcal{L}\left[ \sum_{k=1}^{4}\frac{\exp \left( r_{k}^{2}\,t\right) \,%
\mathrm{erfc}\left( -r_{k}\,\sqrt{t}\right) }{\left( 4a\,r_{k}+3b\right)
r_{k}};s\right] \right\} .  \label{yf_a}
\end{eqnarray}%
Now, apply the convolution theorem of the Laplace transform \cite[Theorem
2.39]{Schiff} to get%
\begin{equation}
y_{f}\left( t\right) =\int_{0}^{t}f\left( t-\tau \right) \sum_{k=1}^{4}\frac{%
\exp \left( r_{k}^{2}\,\tau \right) \,\mathrm{erfc}\left( -r_{k}\,\sqrt{\tau
}\right) }{\left( 4a\,r_{k}+3b\right) r_{k}}d\tau .  \label{L-1[F]_resultado}
\end{equation}

We summarize the above calculations as follows.

\begin{theorem}
\label{Theorem: general_solution}The solution of the Bagley-Torvik equation
\begin{equation}
a\,y^{\prime \prime }\left( t\right) +b\,_{0}D_{t}^{3/2}y\left( t\right)
+c\,y\left( t\right) =f\left( t\right) ,  \label{BT_eqn}
\end{equation}%
with $a,b,c\in
\mathbb{R}
$, and initial conditions $y\left( 0\right) $, $y^{\prime }\left( 0\right) $%
, is given by%
\begin{equation}
y\left( t\right) =y_{f}\left( t\right) +y_{c}\left( t\right) ,
\label{y=yf+yc}
\end{equation}%
where%
\begin{equation}
y_{f}\left( t\right) =\int_{0}^{t}f\left( t-\tau \right) \left(
\sum_{k=1}^{4}\frac{\exp \left( r_{k}^{2}\,\tau \right) \,\mathrm{erfc}%
\left( -r_{k}\sqrt{\tau }\right) }{\left( 4a\,r_{k}+3b\right) \,r_{k}}%
\right) d\tau ,  \label{yf_resultado}
\end{equation}%
takes into account the effect of the external force $f\left( t\right) $, and
\begin{equation}
y_{c}\left( t\right) =\sum_{k=1}^{4}\frac{\left( a\,r_{k}+b\right) \left[
r_{k}^{2}\,\,y\left( 0\right) +y^{\prime }\left( 0\right) \right] }{\left(
4a\,r_{k}+3b\right) r_{k}^{2}}\exp \left( r_{k}^{2}\,t\right) \,\mathrm{erfc}%
\left( -r_{k}\,\sqrt{t}\right) ,  \label{yc_resultado}
\end{equation}%
takes into account the effect of the initial conditions. Also, $r_{k}$ $%
\left( k=1,2,3,4\right) $ are the $4$ different roots of the polynomial $%
P\left( r\right) =a\,r^{4}+b\,r^{3}+c$, which are given by (\ref{Parameters_quartic}) and (\ref{rk_roots}).

\end{theorem}

\begin{remark}
Note that for homogeneous initial conditions, i.e. $y\left( 0\right)
=y^{\prime }\left( 0\right) =0$, the solution of the Bagley-Torvik equation (%
\ref{BT_eqn})\ is reduced to $y\left( t\right) =y_{f}\left( t\right) $ since
$y_{c}\left( t\right) =0$.
\end{remark}

\section{The potential solution \label{Section: potential solution}}

Consider (\ref{BT_equation})\ with $f\left( t\right) =C\,\,t^{\alpha }$,
where $C,\alpha \in
\mathbb{R}
$ are given constants. Note that
\begin{equation}
F\left( s\right) =\mathcal{L}\left[ f\left( t\right) ;s\right] =C\,\mathcal{L%
}\left[ t^{\alpha };s\right] =C\,\frac{\Gamma \left( \alpha +1\right) }{%
s^{\alpha +1}}.  \label{Potential_1}
\end{equation}%
Now, let us calculate (\ref{yf_def}), applying (\ref%
{Inverse_Laplace_transform}) for $\lambda =-1-\alpha $, $p=3$, $q=2$, and $%
m=2$,
\begin{eqnarray}
y_{f}\left( t\right) &=&\mathcal{L}^{-1}\left[ \frac{F\left( s\right) }{%
a\,s^{2}+b\,s^{3/2}+c};t\right]  \notag \\
&=&C\,\Gamma \left( \alpha +1\right) \,\mathcal{L}^{-1}\left[ \frac{%
s^{-1-\alpha }}{a\,s^{2}+b\,s^{3/2}+c};t\right]  \notag \\
&=&C\,\Gamma \left( \alpha +1\right) \,t^{1/2+\alpha }\sum_{k=1}^{4}\frac{%
\mathrm{E}_{1/2,3/2+\alpha }\left( r_{k}\,\sqrt{t}\right) }{%
4a\,r_{k}^{3}+3b\,r_{k}^{2}}.  \label{Potential_2}
\end{eqnarray}

We generalize the above result in the following Theorem.

\begin{theorem}
\label{Theorem: potential solution}The solution of the Bagley-Torvik
equation
\begin{equation}
a\,y^{\prime \prime }\left( t\right) +b\,_{0}D_{t}^{3/2}y\left( t\right)
+c\,y\left( t\right) =\sum_{\ell =0}^{n}C_{\ell }\,t^{\alpha _{\ell }},
\end{equation}%
with $a,b,c,C_{\ell },\alpha _{\ell }\in
\mathbb{R}
$ and initial conditions $y\left( 0\right) $, $y^{\prime }\left( 0\right) $,
is given by
\begin{equation}
y\left( t\right) =y_{f}\left( t\right) +y_{c}\left( t\right) ,
\end{equation}%
where $y_{c}\left( t\right) $ is given by (\ref{yc(t)_resultado})\ and
\begin{equation}
y_{f}\left( t\right) =\sum_{\ell =0}^{n}C_{\ell }\,\Gamma \left( \alpha
_{\ell }+1\right) \,t^{1/2+\alpha _{\ell }}\sum_{k=1}^{4}\frac{\mathrm{E}%
_{1/2,3/2+\alpha _{\ell }}\left( r_{k}\,\sqrt{t}\right) }{%
4a\,r_{k}^{3}+3b\,r_{k}^{2}},  \label{y_f(t)_potential}
\end{equation}%
and where recall that $r_{k}$ $\left( k=1,2,3,4\right) $ are the $4$
different roots of the polynomial $P\left( r\right) =a\,r^{4}+b\,r^{3}+c$,
which are given by (\ref{Parameters_quartic}) and (\ref{rk_roots}).
\end{theorem}

For the particular case of a constant force, i.e. $n=0$ and $\alpha _{0}=0$,
we obtain%
\begin{equation}
y_{f}\left( t\right) =C_{0}\sqrt{t}\sum_{k=1}^{4}\frac{\mathrm{E}%
_{1/2,3/2}\left( r_{k}\,\sqrt{t}\right) }{4a\,r_{k}^{3}+3b\,r_{k}^{2}}.
\label{y(t)_constant_ML}
\end{equation}%
Taking into account (\ref{ML(1/2,3/2)}) 
, rewrite (\ref{y(t)_constant_ML})\ as%
\begin{equation}
y_{f}\left( t\right) =C_{0}\sum_{k=1}^{4}\frac{\mathrm{\exp }\left(
r_{k}^{2}\,t\right) \,\mathrm{erfc}\left( -r_{k}\,\sqrt{t}\right) -1}{%
r_{k}^{3}\left( 4a\,r_{k}+3b\right) }.  \label{yf(t)_constant}
\end{equation}

\section{The sinusoidal solution \label{Section: sinusoidal solution}}

Now, consider (\ref{BT_equation})\ with $f\left( t\right) =\Omega \,\sin
\omega \,t$, where $\Omega ,\omega \in
\mathbb{R}
$, i.e. a sinusoidal external force. According to (\ref{yf_resultado}), we
have%
\begin{equation}
y_{f}\left( t\right) =\Omega \sum_{k=1}^{4}\frac{1}{4a\,r_{k}^{2}+3b\,r_{k}}%
\int_{0}^{t}\sin \omega \left( t-\tau \right) \exp \left( r_{k}^{2}\,\tau
\right) \,\mathrm{erfc}\left( -r_{k}\sqrt{\tau }\right) d\tau .
\label{y(t)_general_sinusoidal}
\end{equation}%
Insert in (\ref{y(t)_general_sinusoidal})\ the integral calculated in Appendix \ref{secA1}, i.e. 
\pagebreak
\begin{eqnarray}
&&\int_{0}^{t}\sin \omega \left( t-\tau \right) \exp \left( r_{k}^{2}\,\tau
\right) \,\mathrm{erfc}\left( -r_{k}\sqrt{\tau }\right) d\tau
\label{Sinusoidal_1} \\
&=&\frac{1}{\omega ^{2}+r_{k}^{4}}\left\{ \omega \left[ \exp \left(
r_{k}^{2}t\right) \,\mathrm{erfc}\left( -r_{k}\sqrt{t}\right) -\cos \omega t%
\right] -\,r_{k}^{2}\sin \omega t\right.  \notag \\
&&+\left. \sqrt{\frac{2}{\omega }}r_{k}\left[ \mathrm{S}\left( \sqrt{\frac{%
2\omega }{\pi }t}\right) \left( r_{k}^{2}\cos \omega t-\omega \sin \omega
t\right) -\mathrm{C}\left( \sqrt{\frac{2\omega }{\pi }t}\right) \left(
r_{k}^{2}\sin \omega t+\omega \cos \omega t\right) \right] \right\} ,  \notag
\end{eqnarray}%
where $\mathrm{S}\left( z\right) $ and $\mathrm{C}\left( z\right) $ denote
the Fresnel integrals (see Appendix \ref{secA1} for the definition of the Fresnel integrals). 
After some algebraic manipulations, we arrive at the following Theorem.

\begin{theorem}
\label{Theorem: sinusoidal solution}The solution of the Bagley-Torvik
equation
\begin{equation}
a\,y^{\prime \prime }\left( t\right) +b\,_{0}D_{t}^{3/2}y\left( t\right)
+c\,y\left( t\right) =\Omega\,\sin \omega \,t,
\end{equation}%
with $a,b,c,\Omega,\omega \in
\mathbb{R}
$ and initial conditions $y\left( 0\right) $, $y^{\prime }\left( 0\right) $,
is given by
\begin{equation}
y\left( t\right) =y_{f}\left( t\right) +y_{c}\left( t\right) ,
\end{equation}%
where $y_{c}\left( t\right) $ is given by (\ref{yc(t)_resultado})\ and\
\begin{eqnarray}
&&y_{f}\left( t\right)   \label{yf(t)_sinusoidal} \\
&=&\Omega\,\left[ \omega \left( \sum_{k=1}^{4}\frac{\exp \left( r_{k}^{2}t\right)
\,\mathrm{erfc}\left( -r_{k}\sqrt{t}\right) }{r_{k}\left(
4a\,r_{k}+3b\right) \left( \omega ^{2}+r_{k}^{4}\right) }-B_{-1}\cos \omega
t\right) -B_{1}\sin \omega t\right]   \notag \\
&&+\Omega\,\sqrt{\frac{2}{\omega }}\left[ \mathrm{S}\left( \sqrt{\frac{2\omega }{%
\pi }t}\right) \left( B_{2}\cos \omega t-B_{0}\,\omega \sin \omega t\right) -%
\mathrm{C}\left( \sqrt{\frac{2\omega }{\pi }t}\right) \left( B_{2}\sin
\omega t+B_{0}\,\omega \cos \omega t\right) \right],  \notag
\end{eqnarray}%
where we have defined the constants
\begin{equation}
B_{m}=\sum_{k=1}^{4}\frac{r_{k}^{m}}{\left( 4a\,r_{k}+3b\right) \left(
\omega ^{2}+r_{k}^{4}\right) },  \label{Bm_def}
\end{equation}%
and where recall that $r_{k}$ $\left( k=1,2,3,4\right) $ are the $4$
different roots of the polynomial $P\left( r\right) =a\,r^{4}+b\,r^{3}+c$,
which are given by (\ref{Parameters_quartic}) and (\ref{rk_roots}).
\end{theorem}

\section{Asymptotic behaviour \label{Section: asymptotic}}

\subsection{Asymptotic behaviour of $y_{c}\left( t\right) $}

Next we obtain the asymptotic behaviour of the solution $y\left( t\right) $
as $t\rightarrow 0^{+}$ and as $t\rightarrow +\infty $ from its Laplace
transform $Y\left( s\right) $ by using the following version of the
Tauberian theorem, (see \cite{Andrade}).

\begin{theorem}
\label{Theorem: Tauberian} Consider that the Laplace transform of a function
$x\left( t\right) $ is given by $X\left( s\right) =\mathcal{L}\left[ x\left(
t\right) ;s\right] $. The asymptotic behaviour of $X\left( s\right) $ as $%
s\rightarrow +\infty $ is given by \
\begin{equation}
X\left( s\right) \approx \mathcal{L}\left[ x_{0}\left( t\right) ;s\right]
,\quad s\rightarrow +\infty ,  \label{Theorem_s->inf}
\end{equation}%
where $x_{0}\left( t\right) $ is the asymptotic behaviour of $x\left(
t\right) $ as $t\rightarrow 0^{+}$. Also, the asymptotic behaviour of $%
X\left( s\right) $ as $s\rightarrow 0^{+}$ is given by \
\begin{equation}
X\left( s\right) \approx \mathcal{L}\left[ x_{\infty }\left( t\right) ;s%
\right] ,\quad s\rightarrow 0^{+},  \label{Theorem_s->0}
\end{equation}%
where $x_{\infty }\left( t\right) $ is the asymptotic behaviour of $x\left(
t\right) $ as $t\rightarrow +\infty $.
\end{theorem}

According to (\ref{yc_def}), we have that the Laplace transform of $%
y_{c}\left( t\right) $ is%
\begin{equation}
Y_{c}\left( s\right) =\frac{y\left( 0\right) \left[ a\,s+b\,s^{1/2}\right]
+y^{\prime }\left( 0\right) \left[ a+b\,s^{-1/2}\right] }{%
a\,s^{2}+b\,s^{3/2}+c}.  \label{Yc(s)}
\end{equation}%
Note that
\begin{eqnarray}
Y_{c}\left( s\right)  &\approx &\frac{y\left( 0\right) }{s}\left( 1-\frac{c}{%
a\,s^{2}}\right) +\frac{y^{\prime }\left( 0\right) }{s^{2}},\quad
s\rightarrow +\infty ,  \label{Yc(s)_s->inf} \\
Y_{c}\left( s\right)  &\approx &\,\frac{y^{\prime }\left( 0\right) }{c}%
\left( \frac{b}{\sqrt{s}}+a\right) +\frac{b}{c}y\left( 0\right) \sqrt{s}%
\,,\quad s\rightarrow 0^{+},  \label{Yc(s)_s->0}
\end{eqnarray}%
thus, applying the inverse Laplace transform \cite[Eqn. 2.2]{Schiff},
\begin{equation}
\mathcal{L}^{-1}\left[ s^{\eta };t\right] =\frac{t^{-1-\eta }}{\Gamma \left(
-\eta \right) },  \label{Inverse_Laplace_s^a}
\end{equation}%
and Theorem \ref{Theorem: Tauberian}, we obtain%
\begin{eqnarray}
y_{c}\left( t\right)  &\approx &y\left( 0\right) \left( 1-\frac{c}{2a}%
t^{2}\right) +y^{\prime }\left( 0\right) \,t,\quad t\rightarrow 0^{+},
\label{yc(t)_t->0} \\
y_{c}\left( t\right)  &\approx &\,\frac{b}{c\sqrt{\pi \,t}}\left( y^{\prime
}\left( 0\right) -\frac{y\left( 0\right) }{2t}\right) ,\quad t\rightarrow
+\infty .  \label{yc(t)_t->inf}
\end{eqnarray}

\subsection{Asymptotic behaviour of $y_{f}\left( t\right) $}

On the one hand, let us assume that for $t\rightarrow 0^{+}$, $f\left(
t\right) $ can be approximated by its truncated Maclaurin series:
\begin{equation}
f\left( t\right) \approx \sum_{n=0}^{N}\frac{f^{\left( n\right) }\left(
0\right) }{n!}t^{n},\quad t\rightarrow 0^{+},  \label{McLaurin}
\end{equation}%
thus (\ref{yf_resultado})\ is approximated as%
\begin{equation}
y_{f}\left( t\right) \approx \sum_{k=1}^{4}\frac{1}{\left(
4a\,r_{k}+3b\right) r_{k}}\sum_{n=0}^{N}\frac{f^{\left( n\right) }\left(
0\right) }{n!}\mathcal{I}_{n}\left( r_{k},t\right) ,  \label{yf_McLaurin}
\end{equation}%
where we define the following integrals:%
\begin{equation}
\mathcal{I}_{n}\left( \nu ,t\right) :=\int_{0}^{t}\left( t-\tau \right)
^{n}\exp \left( \nu ^{2}\,\tau \right) \,\mathrm{erfc}\left( -\nu \,\sqrt{%
\tau }\right) d\tau .  \label{In_def}
\end{equation}%
According to \cite[Eqn. 41:10:7]{Atlas}
\begin{equation}
\mathcal{I}_{0}\left( \nu ,t\right) =\frac{1}{\nu ^{2}}\left[ \exp \left(
\nu ^{2}t\right) \,\mathrm{erfc}\left( -\nu \,\sqrt{t}\right) -\left( 1+%
\frac{2\nu }{\sqrt{\pi }}t^{1/2}\right) \right] ,  \label{I0_resultado}
\end{equation}%
and integrating by parts (\ref{In_def}) for $n=1$, and taking into account (\ref{I0_resultado}), we arrive at%
\begin{equation}
\mathcal{I}_{1}\left( \nu ,t\right) =\frac{1}{\nu ^{4}}\left[ \exp \left(
\nu ^{2}t\right) \,\mathrm{erfc}\left( -\nu \,\sqrt{t}\right) -\left( 1+%
\frac{2\nu }{\sqrt{\pi }}t^{1/2}+\nu ^{2}t+\frac{4\nu ^{3}}{3\sqrt{\pi }}%
t^{3/2}\right) \right] .  \label{I1_resultado}
\end{equation}%
Consequently, the first order approximation is given by%
\begin{equation}
y_{f}\left( t\right) \approx \sum_{k=1}^{4}\frac{f\left( 0\right) \mathcal{I}%
_{0}\left( r_{k},t\right) +f^{\prime }\left( 0\right) \mathcal{I}_{1}\left(
r_{k},t\right) }{\left( 4a\,r_{k}+3b\right) r_{k}},\quad t\rightarrow 0^{+}.
\label{yf(t)_t->0}
\end{equation}%
It is worth noting that when $f\left( t\right) $ is a potential-type force,
i.e. $f\left( t\right) =\sum_{\ell =0}^{n}C_{\ell }t^{\alpha _{\ell }}$, the
approximation given in (\ref{McLaurin})\ cannot be applied in general.
However, we can use the definition of the two-parameter Mittag-Leffler
function (\ref{ML_2parameter_def}) 
in order to see that
\begin{equation}
\mathrm{E}_{\alpha ,\beta }\left( z\right) \approx \sum_{n=0}^{N}\frac{z^{n}%
}{\Gamma \left( \alpha n+\beta \right) },\quad z\rightarrow 0.
\label{ML_t->0}
\end{equation}%
Therefore, the asymptotic behaviour of the potential solution given in
Theorem \ref{Theorem: potential solution} becomes%
\begin{equation}
y_{f}\left( t\right) \approx \sum_{\ell =0}^{n}C_{\ell }\,\Gamma \left(
\alpha _{\ell }+1\right) \,t^{1/2+\alpha _{\ell }}\sum_{n=0}^{N}\frac{%
A_{n-2}\,t^{n/2}}{\Gamma \left( \frac{n}{2}+\frac{3}{2}+\alpha _{\ell
}\right) },\quad t\rightarrow 0^{+},  \label{yf_potential_t->0_general}
\end{equation}%
where the constants $A_{\ell }$ are calculated according to (\ref{A_l_def}).
Nevertheless, from (\ref{A_l=0}), the first non-vanishing term in (%
\ref{yf_potential_t->0_general})\ corresponds to $n=3$, thereby
\begin{equation}
y_{f}\left( t\right) \approx A_{1}\sum_{\ell =0}^{n}\frac{C_{\ell
}\,t^{2+\alpha _{\ell }}}{\left( \alpha _{\ell }+1\right) \left( \alpha
_{\ell }+2\right) },\quad t\rightarrow 0^{+}.  \label{yf_potential_t->0}
\end{equation}

On the other hand, consider the asymptotic behaviour of the two-parameter
Mittag-Leffler function \cite[4.4.17]{GorenfloML}:
\begin{equation}
\mathrm{E}_{\alpha ,\beta }\left( z\right) =-\sum_{n=1}^{N-1}\frac{z^{-n}}{%
\Gamma \left( \beta -\alpha n\right) }+O\left( \left\vert z\right\vert
^{-N}\right) ,\quad |z|\rightarrow \infty ,  \label{ML_t->inf}
\end{equation}%
thus taking the term $n=\ 1$ in (\ref{ML_t->inf}), the asymptotic behaviour
of the potential solution given in Theorem \ref{Theorem: potential solution}
becomes
\begin{equation}
y_{f}\left( t\right) \approx -A_{-3}\sum_{\ell =0}^{n}C_{\ell }\,t^{\alpha
_{\ell }},\quad t\rightarrow +\infty ,  \label{yf_potential_asymp}
\end{equation}%
As a consistency test, we are going to derive in an alternative way the
asymptotic behavior as $t\rightarrow +\infty $ for a special case of the
potential force, i.e. a constant force. For this purpose, consider the
following properties of the complementary error function:%
\begin{eqnarray}
\mathrm{erfc}\left( -x\right)  &=&2-\mathrm{erfc}\left( x\right) ,
\label{erfc_reflection} \\
\mathrm{erfc}\left( x\right)  &\approx &\frac{\exp \left( -x^{2}\right) }{%
\sqrt{\pi }x},\quad x\rightarrow +\infty ,  \label{erfc_x->inf}
\end{eqnarray}%
thus%
\begin{equation}
\exp \left( \alpha ^{2}t\right) \,\mathrm{erfc}\left( -\alpha \,\sqrt{t}%
\right) \approx 2\exp \left( \alpha ^{2}t\right) ,\quad t\rightarrow +\infty
.  \label{exp*erfc_t->inf}
\end{equation}%
Take into account (\ref{exp*erfc_t->inf}) in the solution given in (\ref%
{yf(t)_constant}) for the constant force solution, thereby%
\begin{equation}
y_{f}\left( t\right) \approx C_{0}\left( 2\sum_{k=1}^{4}\frac{\mathrm{\exp }%
\left( r_{k}^{2}\,t\right) \,}{r_{k}^{3}\left( 4a\,r_{k}+3b\right) }%
-A_{-3}\right) ,\quad t\rightarrow +\infty ,  \label{yf_constant_asymp_1}
\end{equation}
Nevertheless, since we know that the solution is non-divergent, the
exponential term in (\ref{yf_constant_asymp_1})\ should vanish, and we have
\begin{equation}
y_{f}\left( t\right) \approx -C_{0}\,A_{-3},\quad t\rightarrow +\infty ,
\label{yf_constant_asymp}
\end{equation}%
which agrees with (\ref{yf_potential_asymp}). Similarly, the asymptotic
behaviour of the sinusoidal solution is given by%
\begin{eqnarray}
&&y_{f}\left( t\right)   \label{yf_sinusoidal_asymp} \\
&\approx &\Omega \sqrt{\frac{2}{\omega }}\left[ \mathrm{S}\left( \sqrt{\frac{%
2\omega }{\pi }t}\right) \left( B_{2}\cos \omega t-B_{0}\,\omega \sin \omega
t\right) -\mathrm{C}\left( \sqrt{\frac{2\omega }{\pi }t}\right) \left(
B_{2}\sin \omega t+B_{0}\,\omega \cos \omega t\right) \right]   \notag \\
&&- \Omega \,\left( B_{-1}\,\omega \cos \omega t+B_{1}\sin \omega t\right)
,\quad t\rightarrow +\infty .  \notag
\end{eqnarray}

\section{Numerical examples \label{Section: Numerical examples}}

In order to compare the solution given in Theorem \ref{Theorem:
general_solution} to other analytical solutions found in the literature, we
have set the following values for the parameters of the Bagley-Torvik
equation (\ref{BT_equation}):%
\begin{equation}
a=1.3,\quad b=2.6,\quad c=3.4.  \label{a,b,c_num}
\end{equation}

Note that the analytical solutions found in the literature, i.e. (\ref%
{yf(t)_Podlubny})-(\ref{yc_Pang}), are quite hard to numerically evaluate.
Indeed, the function $\mathrm{E}_{\lambda ,\mu }^{\left( k\right) }\left(
z\right) $ is given by a series (\ref{Dk_ML_def}), the function $G\left(
t\right) $ is given by another series (\ref{G(t)_Podlubny})\ which involves $%
\mathrm{E}_{\lambda ,\mu }^{\left( k\right) }\left( z\right) $, and then $%
y_{f}\left( t\right) $ is given in terms of a convolution integral which
involves $G\left( t\right) $ and $f\left( t\right) $, i.e. (\ref%
{yf(t)_Podlubny}). Fortunately, the series involved in the computation of $%
G\left( t\right) $ and $\mathrm{E}_{\lambda ,\mu }^{\left( k\right) }\left(
z\right) $ are alternating series, so we can use the acceleration method
described in \cite{Cohen}\ in order to compute them. Note that is
practically impossible to compute these alternating series without any kind
of acceleration method. Nevertheless, although we have implemented a very quick
acceleration method for the numerical evaluation of (\ref{yf(t)_Podlubny})-(\ref%
{yc_Pang}), the computation of the proposed solution presented in Theorem %
\ref{Theorem: general_solution} is much faster and stable. Table \ref%
{TableComputTime}\ presents the computational time ratio $\chi =t_{\mathrm{%
old}}/t_{\mathrm{new}}$ between the computational time of the method given
in (\ref{yf(t)_Podlubny})-(\ref{yc_Pang}) (denoted as $t_{\mathrm{old}}$),
and the method described in Theorem \ref{Theorem: general_solution} for $%
f\left( t\right) =J_{0}\left( t\right) $; (\ref{yc(t)_resultado}) and (\ref%
{yf(t)_constant})\ for $f\left( t\right) =1$; and Theorem \ref{Theorem:
sinusoidal solution} for $\sin \omega t$ with $\omega = 2.5$ (denoted as $t_{%
\mathrm{new}}$).

\begin{center}
\begin{table}[htbp] \centering%
\begin{tabular}{|c|c|}
\hline
$f\left( t\right) $ & $\mathbf{\chi }$ \\ \hline
$J_{0}\left( t\right) $ & $29.40$ \\
$1$ & $881.0$ \\
$\sin \omega t$ & $2513$ \\ \hline
\end{tabular}%
\caption{Computational time ratio of $y(t)$ with initial conditions
$y(0)=y^{\prime}(0)=1$.}\label{TableComputTime}%
\end{table}%
\end{center}

\begin{figure}[htbp]
\centering\includegraphics[width=0.9\textwidth]{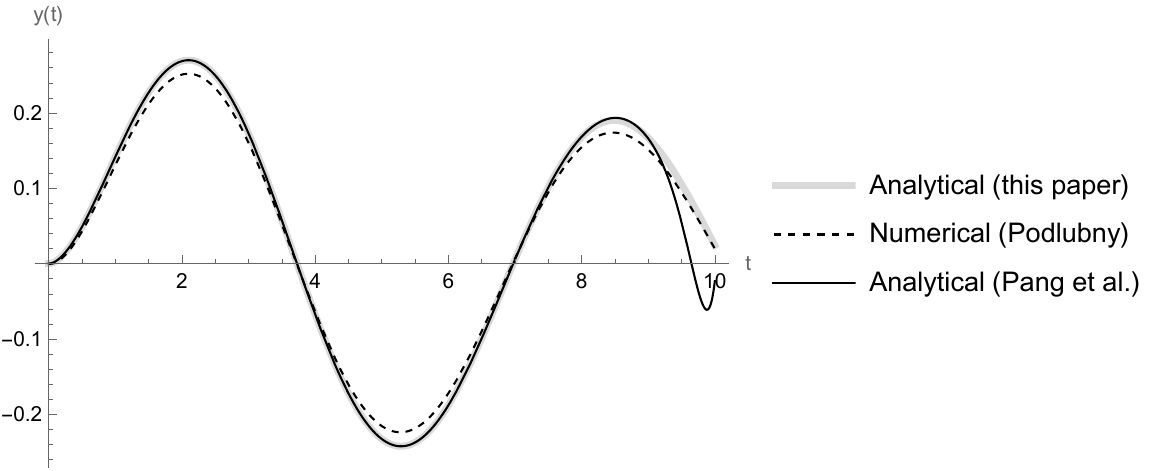}
\caption{Solution of the Bagley-Torvik equation for a Bessel-type force
with $y\left( 0\right) =y^{\prime }\left( 0\right) =0$. }
\label{Figure: BT Bessel homo}
\end{figure}

\begin{figure}[htbp]
\centering\includegraphics[width=0.9\textwidth]{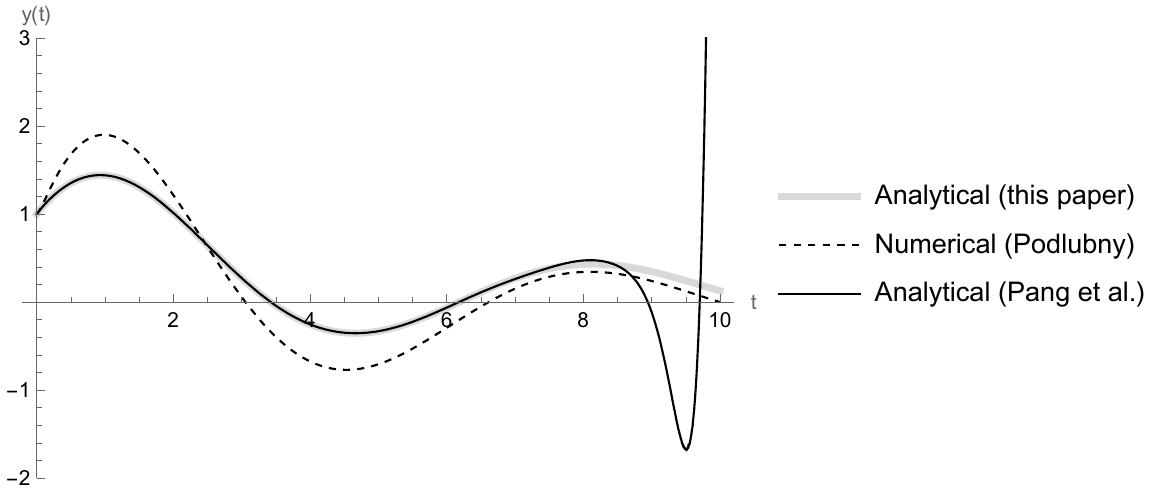}
\caption{Solution of the Bagley-Torvik equation for a Bessel-type force
with $y\left( 0\right) =y^{\prime }\left( 0\right) =1$.}
\label{Figure: BT Bessel nonhomo}
\end{figure}

\begin{figure}[htbp]
\centering\includegraphics[width=0.9\textwidth]{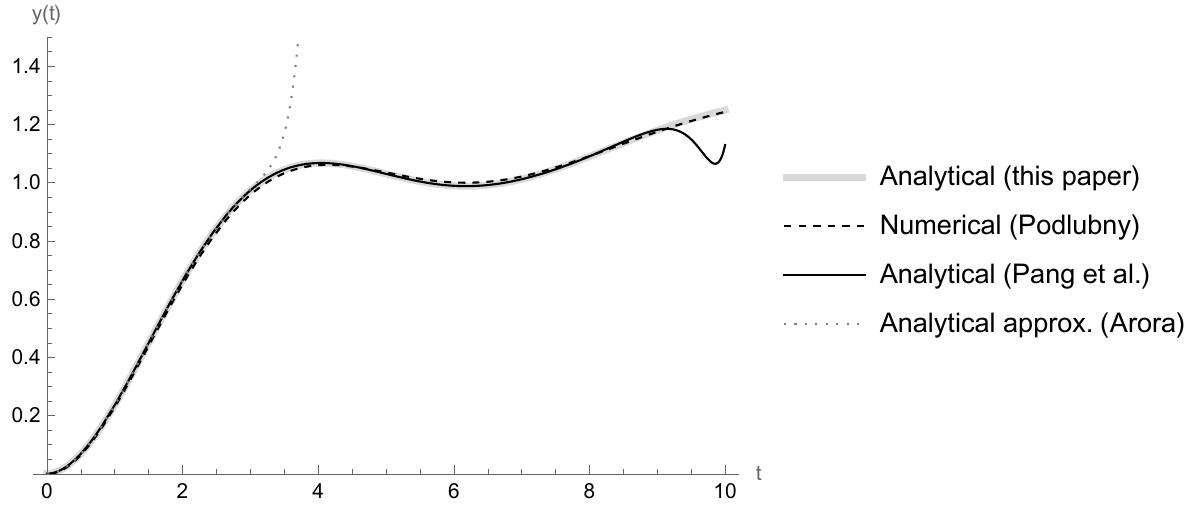}
\caption{Solution of the Bagley-Torvik equation for $f\left( t\right) =1+%
\protect\sqrt{t}$ with $y\left( 0\right) =y^{\prime }\left( 0\right) =0$. }
\label{Figure: BT potential homo}
\end{figure}

\begin{figure}[htbp]
\centering\includegraphics[width=0.9\textwidth]{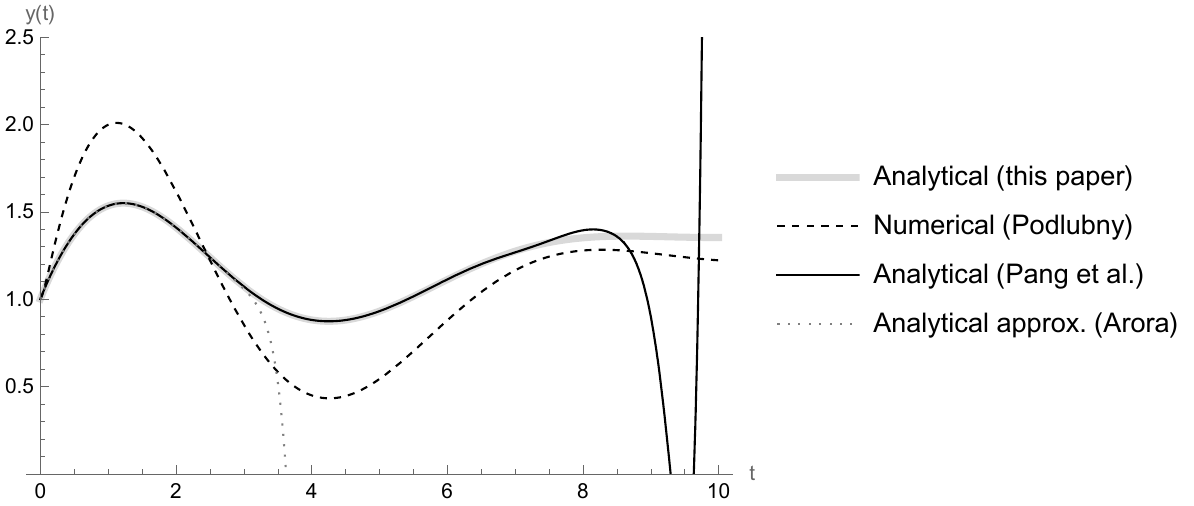}
\caption{Solution of the Bagley-Torvik equation for $f\left( t\right) =1+%
\protect\sqrt{t}$ with $y\left( 0\right) =y^{\prime }\left( 0\right) =1$.}
\label{Figure: BT potential nonhomo}
\end{figure}

\begin{figure}[htb]
\centering\includegraphics[width=0.9\textwidth]{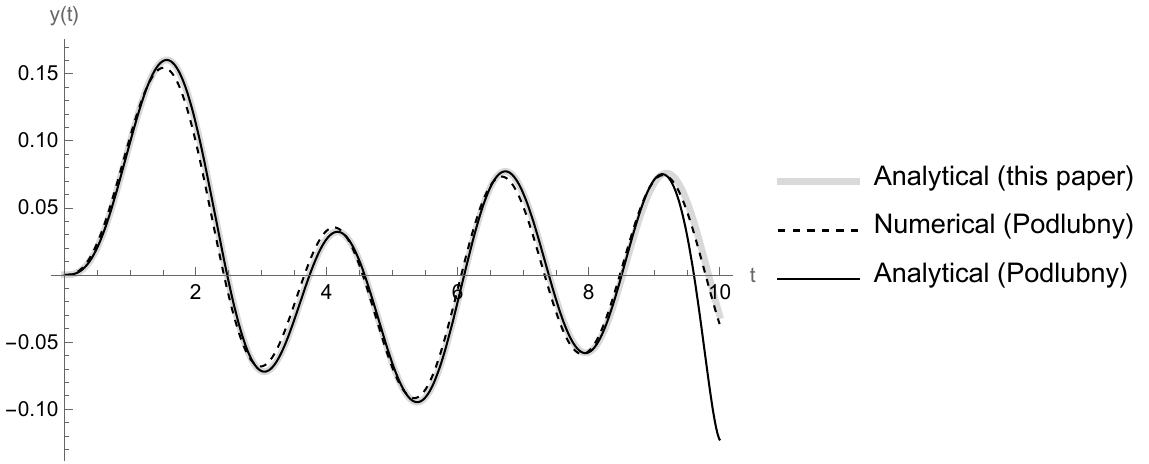}
\caption{Solution of the Bagley-Torvik equation for a sinusoidal force with
$y\left( 0\right) =y^{\prime }\left( 0\right) =0$. }
\label{Figure: BT sinusoidal homo}
\end{figure}

\begin{figure}[htbp]
\centering\includegraphics[width=0.9\textwidth]{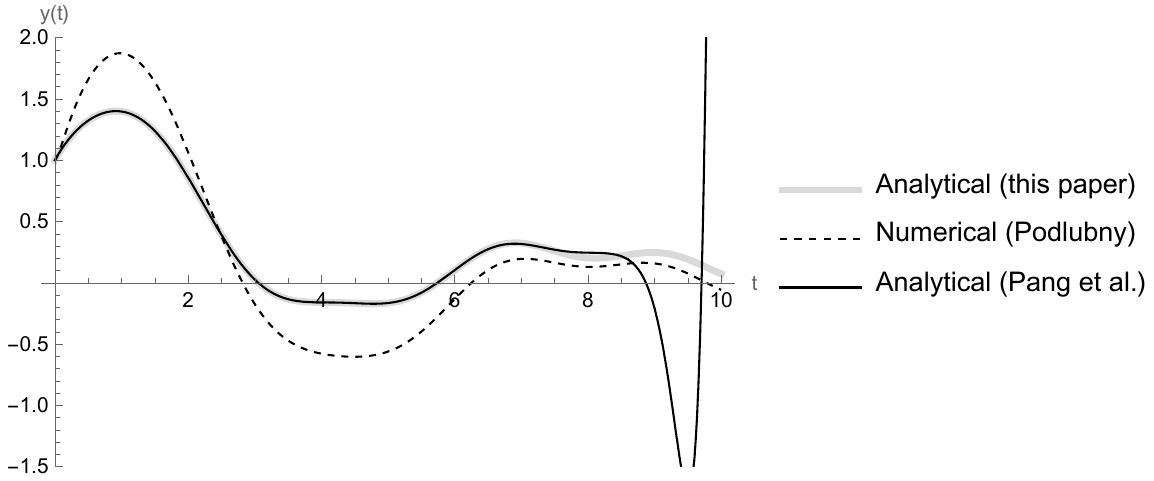}
\caption{Solution of the Bagley-Torvik equation for a sinusoidal force with
$y\left( 0\right) =y^{\prime }\left( 0\right) =1$.}
\label{Figure: BT sinusoidal nonhomo}
\end{figure}

Figures \ref{Figure: BT Bessel homo}-\ref{Figure: BT sinusoidal nonhomo}
compare the solutions of the Bagley-Torvik equation (\ref{BT_equation})\
given in the literature with the solutions proposed in this paper, i.e.
Theorem \ref{Theorem: general_solution} for arbitrary $f\left( t\right) $,
Theorem \ref{Theorem: potential solution}\ for $f\left( t\right) =\sum_{\ell
=0}^{n}C_{\ell }\,t^{\alpha _{\ell }}$, and Theorem \ref{Theorem: sinusoidal
solution}\ for $f\left( t\right) =\Omega \,\sin \omega \,t$. On the one
hand, the analytical solutions found in the literature and computed in this
section are given in (\ref{yf(t)_Podlubny})-(\ref{yc_Pang}), and (\ref%
{f(t)_Arora})-(\ref{Phi(k)_Arora}). The stopping criterion for the recursive
equation (\ref{Recursive_Arora})\ of the latter method is that the series (%
\ref{y(t)_Arora})\ is truncated when $\left\vert d_{k}\right\vert <10^{-12}$%
. On the other hand, the numerical solution found in the literature and
computed in this section is described in Section \ref{Section: Preliminaries}%
, taking $200$ points within the interval $t\in \left( 0,10\right) $.

Figures \ref{Figure: BT Bessel homo} and \ref{Figure: BT
Bessel nonhomo}\ presents the solution of the Bagley-Torvik equation (\ref%
{BT_equation}) with $f\left( t\right) =J_{0}\left( t\right) $, where $%
J_{0}\left( t\right) $ denotes the Bessel function of the first kind of
zeroth order. Also, Figures \ref{Figure: BT potential homo}\ and \ref%
{Figure: BT potential nonhomo} consider $f\left( t\right) =1+\sqrt{t}$, and
Figures \ref{Figure: BT sinusoidal homo}\ and \ref{Figure: BT sinusoidal
nonhomo}\ takes $f\left( t\right) =\sin \omega t$ with $\omega =2.5$.

Figures \ref{Figure: BT Bessel homo}, \ref{Figure: BT
potential homo},\ and \ref{Figure: BT sinusoidal homo}\ consider homogeneous
initial conditions, i.e. $y\left( 0\right) =y^{\prime }\left( 0\right) =0$;
and Figures \ref{Figure: BT Bessel nonhomo}, \ref{Figure: BT potential
nonhomo}\ and \ref{Figure: BT sinusoidal nonhomo}, non-homogeneous initial
conditions, i.e. $y\left( 0\right) =y^{\prime }\left( 0\right) =1$. Note
that the performance of Podlubny's numerical method is reasonably accurate
when we have homogeneous initial conditions, but this is not the case when
we have non-homogeneous initial conditions, where the accuracy is quite low.
Also, Pang's analytical solution is exactly the same as the one proposed in
this paper, except for $t\gtrsim 9$, where the computation of (\ref%
{yf(t)_Podlubny})-(\ref{yc_Pang}) is not stable and starts to oscillate
violently. Similarly, Arora's analytical approximation is exactly the same as
the one proposed here, but for $t\gtrsim 3$, the approximation of the
truncated series (\ref{y(t)_Arora})\ diverges from the actual solution.

\begin{figure}[htbp]
\centering\includegraphics[width=0.9\textwidth]{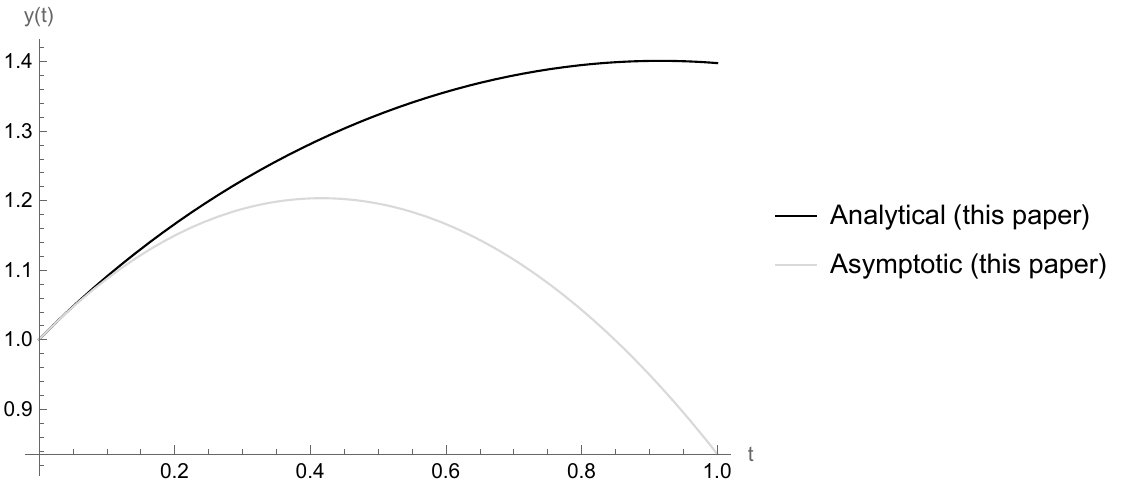}
\caption{Asymptotic behaviour for $t\rightarrow 0^{+}$ of the Bagley-Torvik
equation for a sinusoidal force with $y\left( 0\right) =y^{\prime }\left(
0\right) =1$. }
\label{Figure: Asymp sinusoidal t0}
\end{figure}

\begin{figure}[htbp]
\centering\includegraphics[width=0.9\textwidth]{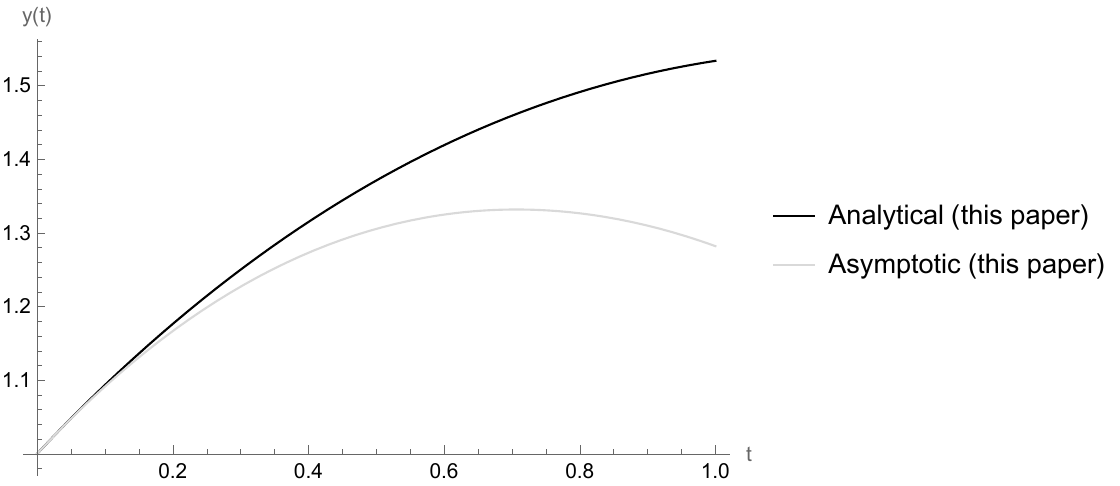}
\caption{Asymptotic behaviour for $t\rightarrow 0^{+}$ of the Bagley-Torvik
equation for $f\left( t\right) =1+\protect\sqrt{t}$ with $y\left( 0\right)
=y^{\prime }\left( 0\right) =1$. }
\label{Figure: Asymp constant t0}
\end{figure}

Figures \ref{Figure: Asymp sinusoidal t0} and \ref{Figure: Asymp sinusoidal
t inf} compares the sinusoidal solution, i.e.
Theorem \ref{Theorem: sinusoidal solution} for $f\left(t\right)=\sin\left(\omega t\right)$ with $\omega = 2.5$,
to its asymptotic behaviour as $t\rightarrow 0^{+}$, i.e. (%
\ref{yc(t)_t->0}) and (\ref{yf(t)_t->0}), and as $t\rightarrow +\infty $, i.e. (%
\ref{yc(t)_t->inf}) and (\ref{yf_sinusoidal_asymp}), respectively. Likewise,
Figures \ref{Figure: Asymp constant t0} and \ref{Figure: Asymp constant t
inf}\ compares the potential solution, i.e. Theorem \ref{Theorem: potential solution} with $f\left( t\right) =1+\sqrt{t}$,
to its asymptotic behaviour as $t\rightarrow 0^{+}$, i.e. (\ref{yc(t)_t->0}) and (\ref%
{yf_potential_t->0}), as well as $t\rightarrow +\infty $, i.e. (\ref{yc(t)_t->inf})
and (\ref{yf_potential_asymp}), respectively. We have taken as initial
conditions $y\left( 0\right) =y^{\prime }\left( 0\right) =1$ in all these
figures.
\begin{figure}[htbp]
\centering\includegraphics[width=0.9\textwidth]{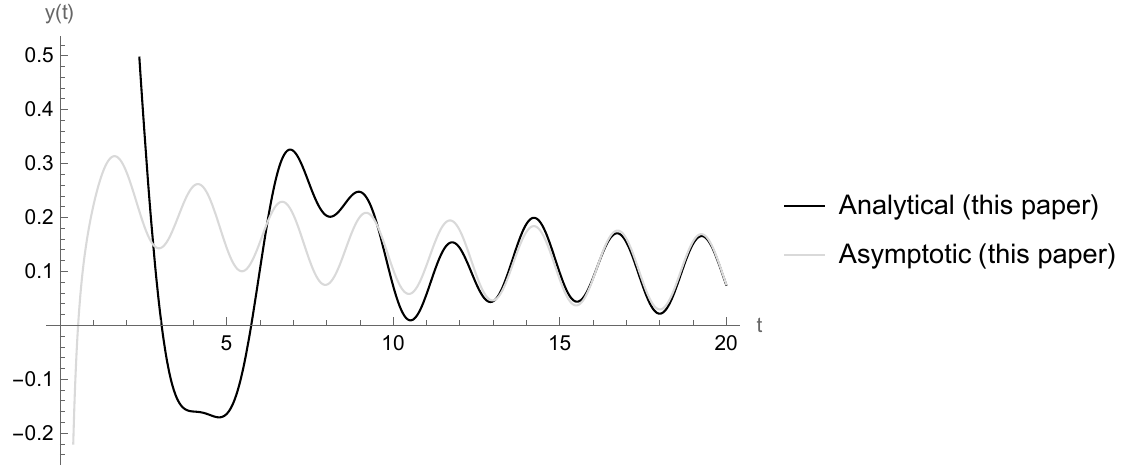}
\caption{Asymptotic behaviour for $t\rightarrow +\infty $ of the
Bagley-Torvik equation for a sinusoidal force with $y\left( 0\right)
=y^{\prime }\left( 0\right) =1$.}
\label{Figure: Asymp sinusoidal t inf}
\end{figure}
\begin{figure}[htbp]
\centering\includegraphics[width=0.9\textwidth]{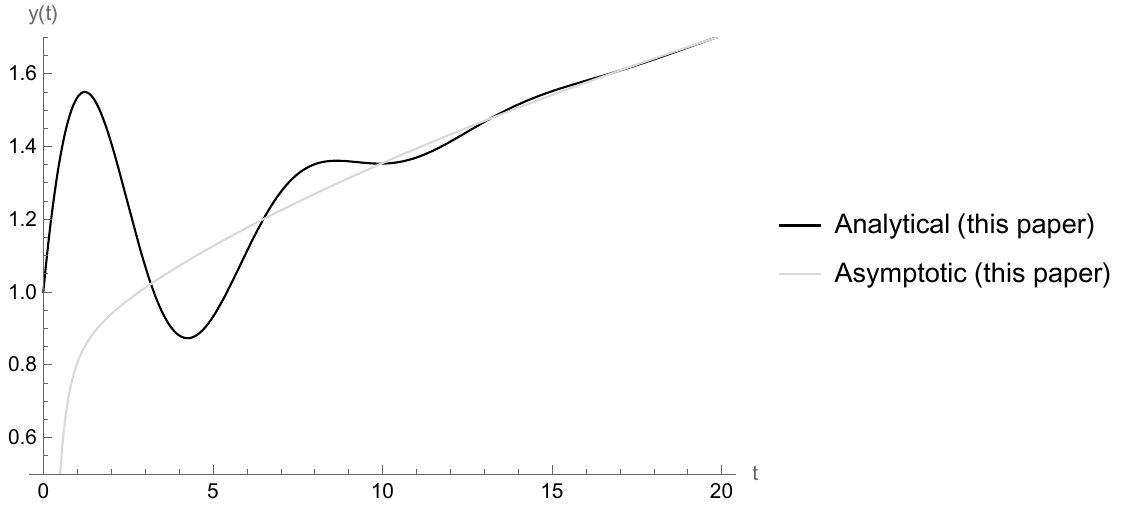}
\caption{Asymptotic behaviour for $t\rightarrow +\infty $ of the
Bagley-Torvik equation for $f\left( t\right) =1+\protect\sqrt{t}$ with $%
y\left( 0\right) =y^{\prime }\left( 0\right) =1$.}
\label{Figure: Asymp constant t inf}
\end{figure}

\pagebreak

\section{Conclusions \label{Section: Conclusions}}

We have calculated the solution $y\left( t\right) $ of the classical
Bagley-Torvik equation (\ref{BT_equation})\ for arbitrary initial conditions
$y\left( 0\right) $, $y^{\prime }\left( 0\right) $, and arbitrary external
force $f\left( t\right) $ by applying the Laplace transform method, and
calculating a new inverse Laplace transform. According to Theorem \ref%
{Theorem: general_solution}, this solution is the sum of two terms $y\left(
t\right) =y_{f}\left( t\right) +y_{c}\left( t\right) $, where $y_{f}\left(
t\right) $\ takes into account the effect of $f\left( t\right) $, and $%
y_{c}\left( t\right) $ takes into account the effect of the initial
conditions. The component $y_{c}\left( t\right) $ is a linear combination of
exponentials with error functions. The component $y_{f}\left( t\right) $ is
expressed in terms of a convolution integral of $f\left( t\right) $ and a
linear combination of exponentials with error functions. The coefficients of
these linear combinations for $y_{c}\left( t\right) $ and $y_{f}\left(
t\right) $ are written in terms of the roots of a quartic equation, whose
analytical expressions are given in (\ref{Parameters_quartic}) and (\ref{rk_roots}).
When $f\left( t\right) $ is a linear combination of
potential or sinusoidal functions, the solution $y\left( t\right) $ can be
expressed in closed-form, i.e. Theorems \ref{Theorem: potential solution}
and \ref{Theorem: sinusoidal solution}, respectively.

On the one hand, we have compared these solutions to other analytical
solutions found in the literature, i.e. (\ref{yf(t)_Podlubny})-(\ref{yc_Pang}%
) and (\ref{f(t)_Arora})-(\ref{Phi(k)_Arora}), and the agreement is
excellent for moderate values of $t$. However, the computation of these
analytical solutions found in the literature is increasingly complicated as $%
t$ increases, and at certain point the numerical accuracy is lost. This is
not the case for the solutions given in Theorems \ref{Theorem:
general_solution}, \ref{Theorem: potential solution}, and \ref{Theorem:
sinusoidal solution}. Also, the numerical computation of the solution given
in these theorems is much faster than the aforementioned analytical
solutions found in the literature, above all in the cases of a sinusoidal or
a potential external force.

On the other hand, we have compared the solutions given in Theorems \ref%
{Theorem: general_solution}, \ref{Theorem: potential solution}, and \ref%
{Theorem: sinusoidal solution} with the numerical solution described in
Section \ref{Section: Preliminaries}, and the agreement is quite good for
homogeneous initial conditions, but is poor for non-homogeneous ones.

In addition, we have calculated the asymptotic behaviour of the solution as $t\rightarrow 0^{+}$
and as $t\rightarrow +\infty$. For the asymptotic behaviour of $y_c\left(t\right)$,
we have used the version of the Tauberian theorem given in Theorem \ref{Theorem: Tauberian}.
In Section \ref{Section: Numerical examples}, we have numerically checked the asymptotic formulas obtained in Section \ref{Section: asymptotic}.

All the computations and plots of this paper have been carried out with
MATHEMATICA. The corresponding MATHEMATICA\ notebook is available at
\hyperref{https://shorturl.at/i0SNN}{}{}{https://shorturl.at/i0SNN}
(accessed on November 2024).

Finally, we hope that the convolution integral given for $y_{f}\left(
t\right) $ can be calculated for other types of external forces $f\left(
t\right) $ as in the case of a sinusoidal force. Also, we expect that the
analytical solution proposed in this paper can be used as a benchmark to
test numerical schemes to solve other types of fractional differential
equations.

\appendix
\section{Auxiliary integrals} \label{secA1}

According to the definition of the Fresnel integrals \cite[Eqns. 7.2.7-8]%
{NIST}%
\begin{eqnarray}
\mathrm{S}\left( z\right) &=&\int_{0}^{z}\sin \left( \frac{\pi t^{2}}{2}%
\right) dt,  \label{Fresnel_S_def} \\
\mathrm{C}\left( z\right) &=&\int_{0}^{z}\cos \left( \frac{\pi t^{2}}{2}%
\right) dt,  \label{Fresnel_C_def}
\end{eqnarray}%
and performing the change of variables $\frac{\pi t^{2}}{2}=\tau $, we obtain%
\begin{eqnarray}
\mathrm{S}\left( \sqrt{\frac{2t}{\pi }}\right) &=&\frac{1}{\sqrt{2\pi }}%
\int_{0}^{t}\frac{\sin \tau }{\sqrt{\tau }}d\tau ,  \label{S_resultado} \\
\mathrm{C}\left( \sqrt{\frac{2t}{\pi }}\right) &=&\frac{1}{\sqrt{2\pi }}%
\int_{0}^{t}\frac{\cos \tau }{\sqrt{\tau }}d\tau .  \label{C_resultado}
\end{eqnarray}
Now, let us calculate the integral:\
\begin{eqnarray}
&&\int_{0}^{t}\sin \left( t-\tau \right) \exp \left( r^{2}\tau \right) \,%
\mathrm{erfc}\left( -r\sqrt{\tau }\right) d\tau  \label{Is_Ic_def} \\
&=&\sin t\underset{I_{c}}{\underbrace{\int_{0}^{t}\cos \tau \exp \left(
r^{2}\tau \right) \,\mathrm{erfc}\left( -r\sqrt{\tau }\right) d\tau }}-\cos t%
\underset{I_{s}}{\underbrace{\int_{0}^{t}\sin \tau \exp \left( r^{2}\tau
\right) \,\mathrm{erfc}\left( -r\sqrt{\tau }\right) d\tau }}.  \notag
\end{eqnarray}%
Integrating by parts, and taking into account (\ref{S_resultado}) and (\ref%
{C_resultado}), we obtain
\begin{eqnarray}
I_{s} &=&\left. -\cos \tau \exp \left( r^{2}\tau \right) \,\mathrm{erfc}%
\left( -r\sqrt{\tau }\right) \right\vert _{0}^{t}+\frac{r}{\sqrt{\pi }}%
\int_{0}^{t}\frac{\cos \tau }{\sqrt{\tau }}d\tau  \label{Is_1} \\
&&+r^{2}\int_{0}^{t}\cos \tau \exp \left( r^{2}\tau \right) \,\mathrm{erfc}%
\left( -r\sqrt{\tau }\right) d\tau  \notag \\
&=&-\cos t\exp \left( r^{2}t\right) \,\mathrm{erfc}\left( -r\sqrt{t}\right)
+1+\sqrt{2}r\ \mathrm{C}\left( \sqrt{\frac{2t}{\pi }}\right)  \notag \\
&&+r^{2}\left[ \sin t\exp \left( r^{2}t\right) \,\mathrm{erfc}\left( -r\sqrt{%
t}\right) -\sqrt{2}r\ \mathrm{S}\left( \sqrt{\frac{2t}{\pi }}\right)
-r^{2}I_{s}\right] ,  \notag
\end{eqnarray}%
thus%
\begin{equation}
I_{s}=\frac{\exp \left( r^{2}t\right) \,\mathrm{erfc}\left( -r\sqrt{t}%
\right) \left( r^{2}\sin t-\cos t\right) +\sqrt{2}r\left[ \mathrm{C}\left(
\sqrt{\frac{2t}{\pi }}\right) -r^{2}\mathrm{S}\left( \sqrt{\frac{2t}{\pi }}%
\right) \right] +1}{1+r^{4}}.  \label{Is_resultado}
\end{equation}%
Similarly, we arrive at%
\begin{equation}
I_{c}=\frac{\exp \left( r^{2}t\right) \,\mathrm{erfc}\left( -r\sqrt{t}%
\right) \left( r^{2}\cos t+\sin t\right) -\sqrt{2}r\left[ \mathrm{S}\left(
\sqrt{\frac{2t}{\pi }}\right) +r^{2}\mathrm{C}\left( \sqrt{\frac{2t}{\pi }}%
\right) \right] -r^{2}}{1+r^{4}}  \label{Ic_resultado}
\end{equation}%
Insert (\ref{Is_resultado})\ and (\ref{Ic_resultado})\ in (\ref{Is_Ic_def}),
and simplify the result
\begin{eqnarray}
&&H_{r}\left( t\right)  \label{H_r(t)} \\
&=&\int_{0}^{t}\sin \left( t-\tau \right) \exp \left( r^{2}\tau \right) \,%
\mathrm{erfc}\left( -r\sqrt{\tau }\right) d\tau  \notag \\
&=&\frac{1}{1+r^{4}}\left\{ \exp \left( r^{2}t\right) \,\mathrm{erfc}\left(
-r\sqrt{t}\right) -r^{2}\sin t-\cos t\right.  \notag \\
&&+\left. \sqrt{2}r\left[ \mathrm{S}\left( \sqrt{\frac{2t}{\pi }}\right)
\left( r^{2}\cos t-\sin t\right) -\mathrm{C}\left( \sqrt{\frac{2t}{\pi }}%
\right) \left( r^{2}\sin t+\cos t\right) \right] \right\} .  \notag
\end{eqnarray}%
Therefore, performing the change of variables $u=\omega \,\tau $, calculate
the following integral as:%
\begin{eqnarray}
&&\int_{0}^{t}\sin \omega \left( t-\tau \right) \exp \left( r^{2}\tau
\right) \,\mathrm{erfc}\left( -r\sqrt{\tau }\right) d\tau
\label{I2_resultado} \\
&=&\frac{1}{\omega }\int_{0}^{\omega t}\sin \left( \omega t-u\right) \exp
\left( \frac{r^{2}}{\omega }u\right) \,\mathrm{erfc}\left( -\frac{r}{\sqrt{%
\omega }}\sqrt{u}\right) du  \notag \\
&=&\frac{1}{\omega }H_{r/\sqrt{\omega }}\left( \omega t\right)  \notag \\
&=&\frac{1}{\omega ^{2}+r^{4}}\left\{ \omega \left[ \exp \left(
r^{2}t\right) \,\mathrm{erfc}\left( -r\sqrt{t}\right) -\cos \omega t\right]
-\,r^{2}\sin \omega t\right.  \notag \\
&&+\left. \sqrt{\frac{2}{\omega }}r\left[ \mathrm{S}\left( \sqrt{\frac{%
2\omega }{\pi }t}\right) \left( r^{2}\cos \omega t-\omega \sin \omega
t\right) -\mathrm{C}\left( \sqrt{\frac{2\omega }{\pi }t}\right) \left(
r^{2}\sin \omega t+\omega \cos \omega t\right) \right] \right\} .  \notag
\end{eqnarray}

\end{document}